\documentclass[a4paper,11pt,reqno]{article}
\usepackage{amsmath, amsfonts, amssymb, amsthm}
\input{epsf.sty}
\usepackage{mathrsfs}
\usepackage{amsmath}
\usepackage{amssymb,color}
\usepackage{graphicx}
\usepackage{subfigure}
\usepackage{amssymb}

\usepackage[colorlinks=true]{hyperref}
\hypersetup{urlcolor=blue, citecolor=blue}

\par

\def \R {{\mathbb{R}}}
\numberwithin{equation}{section}

\begin{document}

\title{A new highly nonlinear equation modelling shallow-water waves\\ with constant vorticity}

\author{Yu $\mbox{Liu}$\footnote{E-mail: TS21080015A31@cumt.edu.cn
}
 \quad Xingxing $\mbox{Liu}$\footnote{Corresponding author. E-mail: liuxxmaths@cumt.edu.cn}\quad
and\quad Min $\mbox{Li}$\footnote{E-mail: TS21080010A31@cumt.edu.cn}\\
$\mbox{School}$ of mathematics, China University of
Mining and Technology,\\Xuzhou, Jiangsu 221116, China}

\date{}
\maketitle

\begin{abstract}
In this paper we apply the approach of formal asymptotic expansions and perturbation theory to derive a new
highly nonlinear shallow-water model from the full governing equations for two dimensional incompressible fluid with constant vorticity. This approximate model is generated by introduction of a
larger scaling than the Camassa-Holm one, which is shown to be optimal in the sense that there are no non-local terms appearing in free surface equation. Moreover, we establish the local well-posedness of the Cauchy problem in Besov spaces, and give a blow-up criterion, which improve the previous corresponding results in Sobolev spaces.\\

\noindent Mathematics Subject Classification: 35Q53, 35G25, 35B30
\smallskip\par
\noindent \textit{Keywords}: Shallow-water waves; Vorticity; Shear flow; Besov spaces; Blow-up.

\end{abstract}

\section{Introduction}
The classical water wave problems embody the Euler equations of motion, the equation of mass conservation and the dynamic and kinematic conditions on the the free surface, and for simplicity the "no-flow" condition on the flat bottom. Due to the complexity of this coupled system of two-dimensional flow with four unknowns (the horizontal and vertical velocity components, the pressure and the free surface), many simpler models have been established as effective approximations in various specific physical regimes. In the derivation of such models, we usually use the following two fundamental dimensionless parameters: the amplitude parameter $\varepsilon=\frac{a}{h_{0}}$ and the shallowness parameter $\mu=\frac{h_{0}^{2}}{\lambda^{2}},$ where $h_0$ is the mean depth of water, $a$ and $\lambda$ are the typical amplitude and wavelength of the waves, respectively. The parameters $\varepsilon$ and $\mu$ clearly represent the contributions of the amplitude and the wavelength of wave under consideration. When speaking about shallow-water waves, we always assume $\mu$ to be small: $\mu\ll 1.$

In the Boussinesq scaling (weakly nonlinear regime): $\mu \ll 1$, $\varepsilon=O(\mu)$,
one can establish the classical weakly nonlinear Korteweg-de Vries (KdV) and Benjamin-Bona-Mahoney (BBM) equations \cite{BBM,KdV}, which are good asymptotic approximations of unidirectional solutions of the two-dimensional irrotational water waves \cite{Constantin-J,Lannes}, but they exhibit behavior that is less nonlinear than dispersive. The two equations play an important role in capturing the existence of solitary water waves \cite{Drazin-J,Johnson1}, and all physically relevant solutions exist globally in time \cite{Constantin0,Tao}. However, the phenomenon of wave breaking (i.e., the solution remains bounded, but its slope becomes unbounded \cite{Whitham}.) is ubiquitous in nature. Thus, in order to gain insight into the most interesting phenomenon of wave breaking, the models derived for the propagation of shallow-water waves require a transition to full nonlinearity.

When allowing $\varepsilon\rightarrow 0$ ($\mu$ fixed), one obtains a set of linear water wave equations. While $\mu\rightarrow0$ ($\varepsilon$ fixed), one finds the fully nonlinear shallow water equations \cite{Constantin-J}. Hence the derivation of the water wave models accommodating wave breaking phenomena amounts to considering larger values of $\varepsilon$. Indeed, the following Camassa-Holm (CH) scaling (moderately nonlinear regime): $\mu \ll 1$, $\varepsilon=O(\sqrt{\mu})$ is proved to be effective. With this scaling, the CH \cite{C-H} and Degasperis-Procesi (DP) equations \cite{DP} are verified as valid approximations to the full governing equations for water waves \cite{Constantin-L}. These two equations are completely integrable Hamiltonian system \cite{Constantin,Constantin-I-L}. More importantly, the CH and DP equations have both solitary waves interacting like solitons \cite{C-H,Constantin-S1,Constantin-S2,L-L-W,L-L,Matsuno}, and in contrast to the KdV and BBM equations, wave breaking solutions \cite{Constantin-E,Constantin-L,Escher-L-Y,L-Y}. Recently, inspired by the idea that the larger $\varepsilon$, the stronger nonlinearity, the models with higher-order nonlinear terms have been derived with larger than the one in the CH scaling. A new highly nonlinear shallow water equation within the scaling of $\mu \ll 1$, $\varepsilon=O(\sqrt[4]{\mu})$, along with the local well-posedness and blow-up criterion were obtained in \cite{Quirchmayr}. Moreover, several integrable equations with cubic nonlinearity as asymptotic models for the propagation of shallow water waves were derived under the scaling of $\mu \ll 1$, $\varepsilon=O(\sqrt[5]{{\mu}^2})$  \cite{C-Hu-L}.

Note that the above obtained equations under the Boussinesq, CH or even larger scalings, are all presented with the assumption of irrotational fluid motion. However, vorticity is of great importance to incorporate the effects of currents and wave-current interactions in the flow \cite{Constantin0}. Nonzero constant vorticity is adapted to tidal flows \cite{T-P} and nonconstant vorticity is the hallmark of highly irregular currents \cite{Constantin-S3}. In the last few years, the problem of water waves with nonzero vorticity and especially with constant vorticity has been given increasing attention \cite{Constantin-E-W,Constantin-E1,Constantin-J2,Constantin-S-S,Constantin-S3,Johnson4,K-W,T-P} and the references therein. Thus the question as to whether these above approximate equations can be appropriate for
the sake of inclusion of an underlying vorticity is an interesting one. R. S. Johnson \cite{Johnson4} showed the CH equation is directly relevant to water waves moving over a
linear shear flow (i.e., constant vorticity in the undisturbed state), and also presented the variation of the depth as the linear shear is varied.
D. Ionescu-Kruse \cite{Ionescu-Kruse} established the CH equation for unidirectional surface shallow water waves moving over a non-zero vorticity flow by an interplay of variational methods and small-parameter expansions.
H. Wang et al. \cite{W-K-L} derived a generalized CH equation with the effect of underlying linear shear flow from the two-dimensional rotational fluid in the CH scaling.
For two component equations modelling shallow-water waves with constant vorticity, we refer the readers to \cite{Escher-H-K-Y,Ivanov}.

In this paper, following \cite{Johnson4}, we formally present the derivation of a highly nonlinear shallow-water model by applying asymptotic expansion method to the
full governing equations for two dimensional flow with an underlying constant vorticity under a larger scaling than the CH one. So a natural question is how to choose an appropriate scaling
that is larger than the CH one and at the same time avoids considerable undertaking. Inspired by the highly nonlinear equations under the assumption of irrotational flow \cite{C-Hu-L,Quirchmayr},
the scalings $\mu \ll 1$, $\varepsilon=O(\sqrt[5]{{\mu}^2})$ and $\mu \ll 1$, $\varepsilon=O(\sqrt[4]{\mu})$ are naturally potential options.
Nevertheless, our calculation confirms that neither of these two scalings is a desired choice. More precisely, on the one hand, within the scaling $\mu \ll 1$, $\varepsilon=O(\sqrt[4]{\mu})$,
it requires us to compute terms up to the order of $O(\varepsilon^4 \mu)$, leaving the $O(\mu^2)$ terms as higher order ones. However, when considering the $O(\varepsilon^4 \mu)$-order approximation,
the term $u_{31,\xi}$ can not be written as a $\xi$-derivative (see Remark \ref{Rem2.1}), which prevents us from getting $u_{31}$ except by means of an integral. Hence we can not obtain the equation
that the surface elevation satisfies in the absence of non-local terms. On the other hand, now that all the terms are available that appear at the order of $O(\varepsilon^3 \mu)$, the scaling  $\mu \ll 1$, $\varepsilon=O(\sqrt[5]{{\mu}^2})$ is certainly usable, but not optimal. Therefore, as mentioned above, we deduce that the scaling $\mu \ll 1$, $\varepsilon=O(\sqrt[3]{\mu})$
is the optimal one larger than the CH in deriving highly nonlinear shallow-water models for two dimensional flow with constant vorticity.

With the above scaling, after deriving the equation of the surface elevation $\eta$, we then establish the following equation of the horizontal velocity $u$ by relating $\eta$ with $u$
{\setlength\arraycolsep{2pt}
\begin{eqnarray}\label{Vorticity}
u_{t}&-&\beta\mu u_{xxt}+cu_{x}+3\alpha\varepsilon uu_{x}-\beta_{0}\mu u_{xxx}+\omega_{1}\varepsilon^{2}u^{2}u_{x}+\omega_{2}\varepsilon^{3}u^{3}u_{x}+\omega_{3}\varepsilon^{4}u^{4}u_{x}+\omega_{4}\varepsilon^{5}u^{5}u_{x}\nonumber\\
&=&\alpha\beta\varepsilon\mu(2u_{x}u_{xx}+uu_{xxx})+\varepsilon^{2}\mu(\omega_{5}u^{2}u_{xxx}+
\omega_{6}u^{3}_{x}+\omega_{7}uu_{x}u_{xx}),
\end{eqnarray}
where the coefficients in Eq. (\ref{Vorticity}) are defined by the right-going wave speed $c=\frac{1}{2}(\sqrt{A^{2}+4}+A)>1,$
with the parameter $A>0$ which represents the vorticity of the underlying flow, $\alpha=\frac{c^4+c^2+1}{3(c^2+1)},\ \beta=\frac{c^4+6c^2+3}{2c^8+6c^6+8c^4+6c^2+2},\ \beta_{0}=\frac{c^6+14c^4+5c^2-2}{6c^9+18c^7+24c^5+18c^3+6c},\
\omega_{1}=\frac{c^5(c^2-1)(c^2+2)}{2(c^2+1)^3},\ \omega_{2}=\frac{c^6(c-1)^2(c+1)^2(c^4+4c^2+6)}{6(c^2+1)^5},\ \omega_{3}=\frac{c^7(c-1)^3(c+1)^3(c^2+4)(c^4+2c^2+6)}{24(c^2+1)^7},\
\omega_{4}=\frac{c^8(c-1)^4(c+1)^4}{120(c^2+1)^9}(c^8+8c^6+28c^4+36c^2+120),\
\omega_{5}=-\frac{2c^{15}+17c^{13}+37c^{11}+115c^9+189c^7+152c^5+54c^3+10c}{12(c^2+c+1)(c^2-c+1)(c^2+1)^5},\
\omega_{6}=\frac{-1}{12(c^2+c+1)(c^2-c+1)(c^2+1)^5}(4c^{16}\\+4c^{15}+26c^{14}+34c^{13}+79c^{12}+137c^{11}+126c^{10}+314c^9+131c^8+402c^7+95c^6+290c^5+50c^4+109c^3+23c^2+18c+6),\
\omega_{7}=\frac{-1}{6(c^2+c+1)(c^2-c+1)(c^2+1)^5}(4c^{16}+2c^{15}+32c^{14}+23c^{13}+115c^{12}+100c^{11}+216c^{10}+319c^9+269c^8+485c^7+233c^6+376c^5+140c^4+137c^3+59c^2+22c+12).$
Applying the transformation $u_{\varepsilon,\mu}(t,x)=\alpha\varepsilon u(\sqrt{\beta\mu}t,\sqrt{\beta\mu}x)$ to Eq. (\ref{Vorticity}), $u_{\varepsilon,\mu}(t,x)$ solves
\begin{eqnarray}\label{equivalent form}
u_{t}&-& u_{txx}+cu_{x}+3uu_{x}-\frac{\beta_{0}}{\beta} u_{xxx}+\frac{\omega_{1}}{\alpha^{2}}u^{2}u_{x}+\frac{\omega_{2}}{\alpha^{3}}u^{3}u_{x}+\frac{\omega_{3}}{\alpha^{4}}u^{4}u_{x}
+\frac{\omega_{4}}{\alpha^{5}}u^{5}u_{x}\nonumber\\
&=&2u_{x}u_{xx}+uu_{xxx}+\frac{1}{\alpha^{2}\beta}(\omega_{7}uu_{x}u_{xx}+\omega_{5}u^{2}u_{xxx}+\omega_{6}u^{3}_{x}).
\end{eqnarray}
Notice that $(1-\partial^2_x)^{-1}f=p\ast f$ for all $f\in L^2(\R)$, where $p(x)=\frac{1}{2}e^{-|x|}$, and $\ast$ denotes convolution with respect to the spatial variable $x$. Then Eq. (\ref{equivalent form}) can be equivalently rewritten as follows
\begin{eqnarray}\label{nonlocal}
u_{t}+(\frac{\beta_{0}}{\beta}+u+\frac{\omega_{5}}{\alpha^2\beta}u^2)u_x&=&p\ast\partial_x\big((\frac{\beta_{0}}{\beta}-c)u-
u^2-\frac{1}{2}u_x^2+\frac{\omega_{5}-\beta \omega_{1} }{3\alpha^2\beta}u^3-\frac{\omega_{2}}{4\alpha^3}u^4-\frac{\omega_{3}}{5\alpha^4}u^5\nonumber\\
&&-\frac{\omega_{4}}{6\alpha^5}u^6+\frac{\omega_{7}-6\omega_{5}}{2\alpha^2\beta}uu^2_{x}\big)+p\ast(\frac{2(\omega_{5}+\omega_{6})-\omega_{7}}{2\alpha^2\beta}u^3_{x}).
\end{eqnarray}
Indeed, we have considered the local well-posedness of the Cauchy problem for Eq. (\ref{nonlocal}) with different coefficients
for initial data $u_0\in H^s(\R),\ s>3/2$ \cite{L-L-L}. Notice that the approaches used in \cite{Danchin,Danchin1} have been successfully applied to the well-posedness of various shallow water wave models in Besov spaces.
Following their idea, our purpose here is to study the local well-posedness in Besov spaces $B_{p,r}^s$ with $1\leq p, r\leq +\infty$ and $s>\max({1+\frac{1}{p}},\frac{3}{2})$.
Thus we recover the well-posedness results in Sobolev spaces $H^s=B_{2,2}^s, \forall s>\frac{3}{2}$ in \cite{L-L-L}. Moreover, we consider the local well-posedness of Eq. (\ref{nonlocal}) for the initial data in $B_{p,r}^s$ with the index $s=\frac{3}{2}$. Noting that $H^s \hookrightarrow B_{2,1}^{\frac{3}{2}}, \forall s>\frac{3}{2}$, we improve our previous local well-posedness results in \cite{L-L-L}. The difficulty we encounter in comparison with CH and DP equations is the complicated higher-order nonlinear structure of Eq. (\ref{nonlocal}), which makes
it difficult to obtain the uniform boundedness of the approximate solutions. To overcome it, we follow the method of dealing with the
Cauchy problem for shallow water waves of large amplitude in \cite{FY}. Then by using commutator estimates in the nonhomogeneous Besov framework, we present a more precise blow-up scenario compared with our previous one given in \cite{L-L-L}.

The remainder of the paper is organized as follows. Section 2 gives the derivation of Eq. (\ref{equivalent form}) from the full governing equations
for two dimensional incompressible fluid with constant vorticity under the scaling $\mu \ll 1$, $\varepsilon=O(\sqrt[3]{\mu})$. Section 3 is devoted to the local well-posedness in Besov spaces.
In Section 4, we derive a new blow-up scenario. In Section 5, we enclose the Littlewood-Paley decomposition, Besov spaces and the transport
equation theory as an Appendix for completeness.

\section{Derivation of the model equation}
\newtheorem{theorem2}{Theorem}[section]
\newtheorem{lemma2}{Lemma}[section]
\newtheorem {remark2}{Remark}[section]
\newtheorem {definition2}{Definition}[section]
\newtheorem{corollary2}{Corollary}[section]
\par
Assume that the two dimensional fluid is incompressible and inviscid with a constant density and zero surface tension. Let $z=0$ denote the location of the flat bottom and $h_0$ be the mean depth, or the undisturbed depth of water. Suppose that $D_{t}=\{(x,z):0<z<h_{0}+\eta(t,x)\}$, where $\eta(t,x)$  measures the deviation from the average level. As in \cite{Johnson1}, the pressure of the fluid is written as $P(t,x,z)=P_{a}+\rho g(h_{0}-z)+p(t,x,z)$,
where $P_{a}$ is the constant atmospheric pressure, the variable $p$ measures the deviation from the hydrostatic pressure distribution, $g$ is the constant
Earth's gravity acceleration and $\rho$ is the constant fluid density. The equations governing
the motion of the fluid consist of Euler equations, together with the equation of mass conservation
\begin{eqnarray*}
u_x+w_z=0, \ \mbox{in} \ D_{t},
\end{eqnarray*}
and the dynamic and kinematic boundary conditions
\begin{align}
\nonumber
\begin{cases}
p=\rho g\eta, &$on$ \ z=h_{0}+\eta(t,x),\\
w=\eta_{t}+u\eta_{x}, &$on$ \ z=h_{0}+\eta(t,x),\\
w=0,  &$on$ \ z=0,
\end{cases}
\end{align}
where $(u(t,x,z),w(t,x,z))$ is the two-dimensional velocity field.

\subsection{Nondimensionalisation and scaling}
Now we introduce the standard dimensionless quantities as in \cite{Constantin-J,Johnson1,Johnson2}, $x\rightarrow\lambda x, z\rightarrow h_{0}z, \eta\rightarrow a\eta, t\rightarrow\frac{\lambda}{\sqrt{gh_{0}}}t,$
and $u\rightarrow\sqrt{gh_{0}}u, w\rightarrow\sqrt{\mu gh_{0}}w, p\rightarrow\rho gh_{0}p.$
With the above nondimensionalised variables, the governing equations become
\begin{align}\label{nondimensionalisation}
\begin{cases}
u_{t}+uu_{x}+wu_{z}=-p_{x}, &$in$ \  0<z<1+\varepsilon\eta(t,x),\\
\mu(w_{t}+uw_{x}+ww_{z})=-p_{z}, &$in$ \ 0<z<1+\varepsilon\eta(t,x), \\
u_{x}+w_{z}=0, &$in$ \ 0<z<1+\varepsilon\eta(t,x),\\
p=\varepsilon\eta, &$on$ \ z=1+\varepsilon\eta(t,x),\\
w=\varepsilon(\eta_{t}+u\eta_{x}), &$on$ \ z=1+\varepsilon\eta(t,x),\\
w=0, &$on$ \ z=0,
\end{cases}
\end{align}
where we have used the two small dimensionless parameters $\varepsilon=a/h_{0}$ and $\mu=h_{0}^{2}/\lambda^{2}.$

Note that $w$ and $p$ (on $z=1+\varepsilon\eta$) are essentially proportional to $\varepsilon$ from the surface boundary conditions in Eq. (\ref{nondimensionalisation}). Then for consistency, we require that $u$ and $w$ are similarly scaled. On the other hand, there is an exact solution of Eq. (\ref{nondimensionalisation}) of the form $(u,w,p,\eta)\equiv (U(z),0,0,0)$, for any $U(z)$. This solution represents laminar flows with a flat free surface and with an arbitrary shear \cite{Johnson4}. In order to consider undulating waves in the presence of a shear flow, we make use of the following scaling
\begin{align}
\nonumber
u\rightarrow U(z)+\varepsilon u, \quad  w\rightarrow\varepsilon w, \quad p\rightarrow\varepsilon p.
\end{align}
Then the governing equations (\ref{nondimensionalisation}) in the presence of a shear flow become
\begin{align}\label{scaling and shear}
\begin{cases}
u_{t}+Uu_x+wU'+\varepsilon(uu_{x}+wu_{z})=-p_{x}, &$in$ \  0<z<1+\varepsilon\eta(t,x),\\
\mu(w_{t}+Uw_x+\varepsilon(uw_{x}+ww_{z}))=-p_{z}, &$in$ \ 0<z<1+\varepsilon\eta(t,x), \\
u_{x}+w_{z}=0, &$in$ \ 0<z<1+\varepsilon\eta(t,x),\\
p=\eta, &$on$ \ z=1+\varepsilon\eta(t,x),\\
w=\eta_{t}+(U+\varepsilon u)\eta_{x}, &$on$ \ z=1+\varepsilon\eta(t,x),\\
w=0, &$on$ \ z=0,
\end{cases}
\end{align}
where the prime in the first equation of (\ref{scaling and shear}) denotes the derivative with respect to $z$. For the case of general $U(z)$, as pointed out in \cite{Johnson4}, the complexity of extensive calculation suggests us to consider a nontrivial special case: linear shear flow. That is, we take $U(z)=Az$, where $A$ is a constant and $0\leq z\leq1$. Here we choose $A>0$, so that the underlying current is in the positive direction of the  $x$-coordinate. The Burns condition \cite{Burns} gives an expression for the speed $c$ of the travelling waves in linear approximation. Specifically, we have
$\int_0^1\frac{dx}{(Az-c)^2}=1,$
which implies that $c=\frac{1}{2}(A\pm\sqrt{A^{2}+4})$. This expression of $c$ have been confirmed again in the leading-order approximation  $O(\varepsilon^{0}\mu^{0})$ \cite{W-K-L}. When $A\equiv0$ (no shear flow), the speed $c=\pm1.$

Before the nondimensionalisation and scaling, the vorticity $\omega=U'+u_z-w_x.$ Using the nondimensionalised variables and the scalings of $u,w$ as before, and the vorticity scaling $\omega\rightarrow\ \sqrt{g/h_0}\omega$, we get $\omega=A+\varepsilon(u_z-\mu w_x)$. Note that the Burns condition arises as a local bifurcation condition for shallow water waves with constant vorticity \cite{Constantin0}. Thus, in order to find a solution with constant vorticity, we require that
\begin{eqnarray}\label{vorticity}
u_z-\mu w_x=0,
\end{eqnarray}
which yields that the vorticity $\omega\equiv A.$

Finally, applying the appropriate far field transformation $\xi=\varepsilon^{\frac{1}{2}}(x-ct),  \tau=\varepsilon^{\frac{3}{2}}t,$
where $c$ is the constant speed for linear propagation, replacing $w$ by $\sqrt{\varepsilon}W$ due to the consistency of the equation of mass conservation,
and then combining the governing equations (\ref{scaling and shear}) with (\ref{vorticity}), we obtain the final version of the governing equations as follows
\begin{equation}\label{final version}
\begin{cases}
\varepsilon u_{\tau}+(Az-c)u_{\xi}+\varepsilon(uu_{\xi}+Wu_{z})+AW=-p_{\xi}, &$in$ \ 0<z<1+\varepsilon\eta(t,x),\\
\varepsilon\mu(\varepsilon W_{\tau}+(Az-c)W_{\xi}+\varepsilon (uW_{\xi}+WW_{z}))=-p_{z}, &$in$ \ 0<z<1+\varepsilon\eta(t,x),\\
u_{\xi}+W_{z}=0, &$in$ \ 0<z<1+\varepsilon\eta(t,x),\\
u_{z}-\varepsilon\mu W_{\xi}=0, &$in$ \ 0<z<1+\varepsilon\eta(t,x),\\
p=\eta, &$on$ \ z=1+\varepsilon\eta(t,x),\\
W=\varepsilon\eta_{\tau}+(A-c)\eta_{\xi}+\varepsilon\eta_{\xi}( u+A\eta), &$on$ \ z=1+\varepsilon\eta(t,x),\\
W=0, &$on$ \ z=0.\\
\end{cases}
\end{equation}

\subsection{A new highly nonlinear equation}
In search for a solution of Eq. (\ref{final version}) formally, we suppose that the concerning functions can be
written as double asymptotic expansions in the two parameters $\varepsilon$ and $\mu$: $f\sim\sum\limits_{i=0}^{\infty}\sum\limits_{j=0}^{\infty}\varepsilon^{i}\mu^{j}f_{ij},$
as $\varepsilon\rightarrow 0$ and $\mu\rightarrow0$, where $f$ (and correspondingly $f_{ij}$) stands for each of $u,W,p,\eta,$ and all the functions $f_{ij}$ satisfy the far field conditions $f_{ij}\rightarrow0$ for every $i,j=0,1,2,3,...$  as $|\xi|\rightarrow\infty$. Note that the relations of $u,W,p$ and $\eta$ need to be satisfied on the free surface $1+\varepsilon\eta$, but itself is unknown. We deal with this difficulty by taking advantage of Taylor expansions of $u,W,p$ on the surface about $z=1$.

Now we can perform formal calculation by substituting the asymptotic expansions of $u,W,p,\eta$ into Eq. (\ref{final version}). We need to compute
all the coefficients of the orders $O(\varepsilon^{i}\mu^{j})(i.j=0,1,2,3,...)$. Indeed, the orders $O(\varepsilon^{0}\mu^{0}),O(\varepsilon^{1}\mu^{0}),O(\varepsilon^{2}\mu^{0}), O(\varepsilon^{3}\mu^{0}),O(\varepsilon^{4}\mu^{0}),O(\varepsilon^{0}\mu^{1}),O(\varepsilon^{1}\mu^{1})$ and $O(\varepsilon^{2}\mu^{1})$ have been given in \cite{W-K-L}.
Thus it remains us to calculate the asymptotic expansions at orders $O(\varepsilon^{5}\mu^{0}),O(\varepsilon^{6}\mu^{0})$ $ \mbox{and}\ O(\varepsilon^{3}\mu^{1}).$
After tedious and lengthy calculations, one can get the desired results as follows
 \begin{align}\label{eta40tau1}
(2c&-A)\eta_{40,\tau}+(3c^2-3Ac+A^2)(\eta_{40}\eta_{00}+\eta_{30}\eta_{10}+\frac{1}{2}\eta_{20}^2)_{\xi}-3(2c+\frac{4}{3}c_1-A)(c+c_1-\frac{A}{2})\nonumber\\
&\cdot(\eta_{30}\eta_{00}^2+2\eta_{20}\eta_{10}\eta_{00}+\frac{1}{3}\eta_{10}^3)_{\xi}+4((2c+\frac{3}{2}c_1-A)(c+c_1-c_2-\frac{A}{2})+\frac{1}{2}(c+c_1-\frac{A}{2})\nonumber\\
&\cdot(c+c_1-3c_2-\frac{A}{2}))(\eta_{00}^3\eta_{20}+\frac{3}{2}\eta_{00}^2\eta_{10}^2)_{\xi}-5((2c+\frac{8}{5}c_1-A)(c+c_1-c_2+c_3-\frac{A}{2})\nonumber\\
&+(c+c_1-\frac{9 c_2}{5}-\frac{A}{2})+\frac{8 c_3}{5}(c+c_1-\frac{A}{2}))(\eta_{00}^4\eta_{10})_{\xi}-((c+c_1-c_2+c_3-c_4-\frac{A}{2})\nonumber\\
&\cdot(2c+\frac{5 c_1}{3}-A)+(c+c_1-2c_2-\frac{A}{2})(c+c_1-c_2+c_3-\frac{A}{2})+\frac{1}{2}(c+c_1-c_2+4c_3-\frac{A}{2})\nonumber\\
&\cdot(c+c_1-c_2-\frac{A}{2})-\frac{5}{3}c_4(c+c_1-\frac{A}{2}))(\eta_{00}^6)_{\xi}=0,
\end{align}
\begin{align}\label{eta21tau}
&(2c-A)\eta_{21,\tau}+(3c^{2}-3Ac+A^{2})(\eta_{21}\eta_{00}+\eta_{20}\eta_{01}+\eta_{11}\eta_{10})_{\xi}-3(2c+\frac{4 c_{1}}{3}-A)(c+c_{1}-\frac{A}{2})\nonumber\\
&(\eta_{11}\eta_{00}^{2}+2\eta_{10}\eta_{01}\eta_{00})_{\xi}+4((2c+\frac{3 c_{1}}{2}-A)(c+c_{1}-c_{2}-\frac{A}{2})+(\frac{c}{2}
+\frac{c_{1}}{2}-\frac{3 c_{2}}{2}-\frac{A}{4})(c+c_{1}\nonumber\\
&-\frac{A}{2}))(\eta_{01}\eta_{00}^{3})_{\xi}+\frac{(c-A)^{2}}{3}\eta_{20,\xi\xi\xi}+(\frac{-c^{2}-8cc_{1}+8c_{1}A+A^{2}}{3}+\frac{(c-A)^{2}(2c_{1}-2c+A)}{3(2c-A)})\nonumber\\
&\cdot(\eta_{10,\xi}\eta_{00,\xi})_{\xi}+(\frac{-c^{2}-4cc_{1}+4Ac_{1}+A^{2}}{3}+\frac{(c-A)^{2}(4c+4c_{1}-2A)}{3(2c-A)})(\eta_{00}\eta_{10,\xi\xi}+\eta_{10}\eta_{00,\xi\xi})_{\xi}\nonumber\\
&+c_{8}\eta_{00,\xi}^{3}+c_{9}\eta_{00}^{2}\eta_{00,\xi\xi\xi}+c_{10}\eta_{00}\eta_{00,\xi}\eta_{00,\xi\xi}=0,
\end{align}
where\\
$c_{8}=\frac{2c^{10}+13c^8+19c^6+38c^4+33c^2+9}{6c^2(c^2+1)^4},c_{9}=\frac{c^{10}+3c^8+2c^6+28c^4+21c^2+5}{6c^2(c^2+1)^4},c_{10}=\frac{3c^{10}+15c^8+13c^6+52c^4+44c^2+11}{3c^2(c^2+1)^4},$
and
 \begin{align}\label{eta50tau}
&(2c-A)\eta_{50,\tau}+(3c^{2}-3Ac+A^{2})(\eta_{00}\eta_{50}+\eta_{10}\eta_{40}+\eta_{20}\eta_{30})_{\xi}-3(2c+\frac{4}{3}c_{1}-A)(c+c_{1}-\frac{A}{2})\nonumber\\
&\cdot(\eta_{00}^{2}\eta_{40}+\eta_{20}^{2}\eta_{00}+2\eta_{30}\eta_{10}\eta_{00})_{\xi}+4((2c+\frac{3}{2}c_{1}-A)(c+c_{1}-c_{2}-\frac{A}{2})+(\frac{c}{2}+\frac{c_{1}}{2}-\frac{3 c_{2}}{2}-\frac{A}{4})\nonumber\\
&\cdot(c+c_{1}-\frac{A}{2}))(\eta_{00}^{3}\eta_{30}+\eta_{10}^{3}\eta_{00}+3\eta_{20}\eta_{10}\eta_{00}^{2})_{\xi}-5((2c+\frac{8}{5}c_{1}-A)(c+c_{1}-c_{2}+c_{3}-\frac{A}{2})\nonumber\\
&+(c+c_{1}-\frac{9}{5}c_{2}-\frac{A}{2})(c+c_{1}-c_{2}-\frac{A}{2})+\frac{8}{5}c_{3}(c+c_{1}-\frac{A}{2}))(\eta_{00}^{4}\eta_{20}+2\eta_{00}^{3}\eta_{10}^{2})_{\xi}+6((2c-A\nonumber\\
&+\frac{5}{3}c_{1})(c+c_{1}-c_{2}+c_{3}-c_{4}-\frac{A}{2})+(c+c_{1}-c_{2}+c_{3}-\frac{A}{2})(c+c_{1}-2c_{2}-\frac{A}{2})+\frac{1}{2}(c+c_{1}\nonumber\\
&-c_{2}-\frac{A}{2})(c+c_{1}-c_{2}+4c_{3}-\frac{A}{2})-\frac{5}{3}c_{4}(c+c_{1}-\frac{A}{2}))(\eta_{10}\eta_{00}^{5})_{\xi}-((2c-A+\frac{12}{7}c_{1})(c+c_{1}\nonumber\\
&-c_{2}+c_{3}-c_{4}+c_{5}-\frac{A}{2}))+(c+c_{1}-\frac{15}{7}c_{2}-\frac{A}{2})(c+c_{1}-c_{2}+c_{3}-c_{4}-\frac{A}{2})-\frac{15}{7}c_{4}(c+c_{1}\nonumber\\
&-c_2-\frac{A}{2})+(c+c_{1}-c_{2}+\frac{16}{7}c_{3}-\frac{A}{2})(c+c_{1}-c_{2}+c_{3}-\frac{A}{2})+\frac{12}{7}c_{5}(c+c_{1}-\frac{A}{2}))(\eta_{00}^{7})_{\xi}=0.
\end{align}
According to the double asymptotic expansions of $\eta$, we take $\eta:=\eta_{00}+\varepsilon\eta_{10}+\varepsilon^{2}\eta_{20}+\varepsilon^{3}\eta_{30}+\varepsilon^{4}\eta_{40}+\varepsilon^{5}\eta_{50}+\mu\eta_{01}
+\varepsilon\mu\eta_{11}+\varepsilon^{2}\mu\eta_{21}+O(\varepsilon^{6},\mu^{2},\varepsilon^{3}\mu).$ Multiplying each of the equations that $\eta_{00},\eta_{10},\eta_{20},\eta_{30},\eta_{01},\eta_{11}$ satisfy, and
(\ref{eta40tau1})-(\ref{eta50tau}) by
$1,\varepsilon,\varepsilon^{2},\varepsilon^{3},\mu,\varepsilon\mu,\varepsilon^{4},\varepsilon^{2}\mu,\varepsilon^{5},$ respectively, and then summating these results, we get the equation of $\eta$ up to the order $O(\varepsilon^{6},\mu^{2},\varepsilon^{3}\mu)$
\begin{multline}\label{guanyuaita}
(2c-A)\eta_{\tau}+(3c^{2}-3Ac+A^{2})\eta\eta_{\xi}+\frac{(c-A)^{2}}{3}\mu\eta_{\xi\xi\xi}+\varepsilon A_{1}\eta^{2}\eta_{\xi}+\varepsilon^{2} A_{2}\eta^{3}\eta_{\xi}+\varepsilon^{3} A_{3}\eta^{4}\eta_{\xi}+\varepsilon^{4} A_{4}\eta^{5}\eta_{\xi}\\ \quad+\varepsilon^{5} A_{5}\eta^{6}\eta_{\xi}=\varepsilon\mu(A_{6}\eta_{\xi}\eta_{\xi\xi}+A_{7}\eta\eta_{\xi\xi\xi})+\varepsilon^{2}\mu(A_{8}\eta_{\xi}^{3}+A_{9}\eta^{2}\eta_{\xi\xi\xi}+A_{10}\eta\eta_{\xi}\eta_{\xi\xi})
+O(\varepsilon^{6},\varepsilon^{3}\mu,\mu^{2}),
\end{multline}
where $A_1=-\frac{c^4+4c^2+1}{2(c^2+1)^2},
A_2=\frac{c^8+6c^6+4c^4+6c^2+1}{3(c^2+1)^4},
A_3=-\frac{c^{12}+8c^{10}+9c^8+24c^6+9c^4+8c^2+1}{4(c^2+1)^6},
A_4=\frac{1}{5(c^2+1)^8}\cdot\\(c^{16}+10c^{14}+16c^{12}+60c^{10}+36c^8+60c^6+16c^4+10c^2+1),
A_5=\frac{-1}{6(c^2+1)^{10}}(c^{20}+12c^{18}+25c^{16}+120c^{14}+100c^{12}+240c^{10}+100c^8+120c^6+25c^4+12c^2+1),
A_6=-\frac{2c^6+4c^4+11c^2+6}{3c^2(c^2+1)^2},A_7=-\frac{c^4+6c^2+3}{3c^2(c^2+1)^2},A_8=-c_{8},
A_9=-c_{9},A_{10}=-c_{10}$.

On the other hand, substituting $\eta_{ij}$ represented by $u_{ij}$ into the double asymptotic expansion (\ref{guanyuaita}), and noticing that
$u:= u_{00}+\varepsilon u_{10}+\varepsilon^{2}u_{20}+\varepsilon^{3}u_{30}+\varepsilon^{4}u_{40}+\varepsilon^{5}u_{50}+\mu u_{01}+\varepsilon\mu u_{11}+\varepsilon^{2}\mu u_{21}+O(\varepsilon^{6},\mu^{2},\varepsilon^{3}\mu)$, we obtain
\begin{eqnarray}\label{aita and u guanxishi}
\eta&=&\frac{1}{c-A}u+\gamma_{1}\varepsilon u^{2}+\gamma_{2}\varepsilon^{2}u^{3}+\gamma_{3}\varepsilon^{3}u^{4}+\gamma_{4}\varepsilon^{4}u^{5}+\gamma_{5}\varepsilon^{5}u^{6}+\gamma_{6}\varepsilon\mu u_{\xi\xi}+\gamma_{7}\varepsilon^{2}\mu u_{\xi}^{2}\nonumber\\
&&+\gamma_{8}\varepsilon^{2}\mu uu_{\xi\xi}+O(\varepsilon^{6},\varepsilon^{3}\mu,\mu^{2}).
\end{eqnarray}
In what follows, we shall derive Eq. (\ref{Vorticity}) of $u$ from the surface equation (\ref{guanyuaita}) of $\eta$.
Thanks to (\ref{guanyuaita})-(\ref{aita and u guanxishi}), a direct calculation gives rise to
\begin{align}\label{guanyu u}
u_{\tau}&+2(c-A)\gamma_{1}\varepsilon uu_{\tau}+3(c-A)\gamma_{2}\varepsilon^{2}u^{2}u_{\tau}+4(c-A)\gamma_{3}\varepsilon^{3}u^{3}u_{\tau}+5(c-A)\gamma_{4}\varepsilon^{4}u^{4}u_{\tau}+6(c-A)\nonumber\\
&\times\gamma_{5}\varepsilon^{5}u^{5}u_{\tau}+(c-A)\gamma_{6}\varepsilon\mu u_{\tau\xi\xi}+2(c-A)\gamma_{7}\varepsilon^{2}\mu u_{\xi}u_{\tau\xi}+(c-A)\gamma_{8}\varepsilon^{2}\mu(u_{\tau}u_{\xi\xi}+uu_{\tau\xi\xi})\nonumber\\
&+\frac{3c^{2}-3Ac+A^{2}}{(2c-A)(c-A)}uu_{\xi}+\frac{(c-A)B_{1}}{2c-A}\varepsilon u^{2}u_{\xi}+\frac{(c-A)B_{2}}{2c-A}\varepsilon^{2}u^{3}u_{\xi}+\frac{(c-A)B_{3}}{2c-A}\varepsilon^{3}u^{4}u_{\xi}\nonumber\\
&+\frac{(c-A)B_{4}}{2c-A}\varepsilon^{4}u^{5}u_{\xi}+\frac{(c-A)B_{5}}{2c-A}\varepsilon^{5}u^{6}u_{\xi}+\frac{(c-A)B_{6}}{2c-A}\varepsilon\mu u_{\xi}u_{\xi\xi}+\frac{(c-A)B_{7}}{2c-A}\varepsilon\mu uu_{\xi\xi\xi}\nonumber\\
&+\varepsilon^{2}\mu(\frac{(c-A)B_{8}}{2c-A}u^{2}u_{\xi\xi\xi}+\frac{(c-A)B_{9}}{2c-A}u_{\xi}^{3}+\frac{(c-A)B_{10}}{2c-A}uu_{\xi}u_{\xi\xi})+\frac{(c-A)^{2}}{3(2c-A)}\mu u_{\xi\xi\xi}\nonumber\\
&=O(\varepsilon^{6},\varepsilon^{3}\mu,\mu^{2}),
\end{align}
where $B_1=\frac{c^3(3c^6+5c^4+2c^2+2)}{2(c^2+1)^2},
B_2=\frac{c^4(c+1)(c-1)(7c^8+24c^6+24c^4+9c^2+6)}{6(c^2+1)^4},
B_3=\frac{c^5(c+1)^2(c-1)^2}{24(c^2+1)^6}(15c^{10}+79c^8+160c^6+142c^4+40c^2+24),
B_4=\frac{c^6(c-1)^3(c+1)^3(31c^{12}+222c^{10}+660c^8+1043c^6+854c^4+190c^2+120)}{120(c^2+1)^8},
B_5=\frac{c^7(c-1)^4(c+1)^4(63c^{14}+573c^{12}+2266c^{10}+5006c^8+6892^6+5748c^4+936c^2+720)}{720(c^2+1)^{10}},
B_6=\frac{(c^2+c+1)(c^2-c+1)}{2}z^2-\frac{1}{6(c^2+1)^2}\cdot\\(c^8+c^6-6c^4-21c^2-9),
B_7=\frac{(c^2+c+1)(c^2-c+1)}{2}z^2-\frac{(c-1)(c+1)(c^2+3)(c^4+3c^2+1)}{6(c^2+1)},
B_8=\frac{z^2}{12(c^2+1)^4} (15c^{13}-6c^{12}+57c^{11}-24c^{10}+87c^9-42c^8+75c^7-42c^6+42c^5-24c^4+12c^3-6c^2)
-\frac{1}{12(c^2+1)^4}(5c^{13}+2C^{12}+27c^{11}+43c^9+46c^8+5c^7+82^6-74c^5+84c^4-44c^3+50c^2-10c+12),
B_9=\frac{c^2(c-1)(c^2+c+1)(c^2-c+1)z^2}{2(c^2+1)}-\frac{c^{13}+c^{12}+3c^{10}-14c^9+8c^8-38c^7+14c^6-60c^5+15c^4-38c^3+10c^2-9c+3}{6(c^2+1)^4},
B_{10}=\frac{c^2z^2}{2(c^2+1)^2}(7c^7-4c^6+13c^5-8c^4+10c^3-8c^2+6c-4)-\frac{7c^{13}+4c^{12}+27c^{11}+18c^{10}+17c^9+62c^8-55c^7+110c^6-166c^5+114c^4-100c^3+70c^2-22c+18}{6(c^2+1)^4}.$
It thus follows from (\ref{guanyu u}) that
\begin{align*}
&\varepsilon^{2}u^{2}u_{\tau}=-\varepsilon^{2}u^{2}(2(c-A)\gamma_{1}\varepsilon uu_{\tau}+3(c-A)\gamma_{2}\varepsilon^{2}u^{2}u_{\tau}+4(c-A)\gamma_{3}\varepsilon^{3}u^{3}u_{\tau}+\frac{3c^{2}-3Ac+A^{2}}{(2c-A)(c-A)}uu_{\xi}\\
&+\frac{(c-A)B_{1}}{2c-A}\varepsilon u^{2}u_{\xi}+\frac{(c-A)B_{2}}{2c-A}\varepsilon^{2}u^{3}u_{\xi}+\frac{(c-A)B_{3}}{2c-A}\varepsilon^{3}u^{4}u_{\xi}+\frac{(c-A)^{2}}{3(2c-A)}\mu u_{\xi\xi\xi})+O(\varepsilon^{6},\varepsilon^{3}\mu,\mu^{2}),
\end{align*}
which implies that
\begin{align*}
&\varepsilon^{2}u^{2}(1+2(c-A)\gamma_{1}\varepsilon u+3(c-A)\gamma_{2}\varepsilon^{2}u^{2}+4(c-A)\gamma_{3}\varepsilon^{3}u^{3})u_{\tau}=-\varepsilon^{2}u^{2}(\frac{3c^{2}-3Ac+A^{2}}{(2c-A)(c-A)}uu_{\xi}\\
&+\frac{(c-A)B_{1}}{2c-A}\varepsilon u^{2}u_{\xi}+\frac{(c-A)B_{2}}{2c-A}\varepsilon^{2}u^{3}u_{\xi}+\frac{(c-A)B_{3}}{2c-A}\varepsilon^{3}u^{4}u_{\xi}+\frac{(c-A)^{2}}{3(2c-A)}\mu u_{\xi\xi\xi})+O(\varepsilon^{6},\varepsilon^{3}\mu,\mu^{2}).
\end{align*}
It is thereby adduced that
{\setlength\arraycolsep{2pt}
\begin{eqnarray*}
\varepsilon^{2}u^{2}u_{\tau}&=&-\varepsilon^{2}u^{2}(1-(2(c-A)\gamma_{1}\varepsilon u+3(c-A)\gamma_{2}\varepsilon^{2}u^{2}+4(c-A)\gamma_{3}\varepsilon^{3}u^{3})+4(c-A)^{2}\gamma_{1}^{2}\varepsilon^{2}u^{2}\nonumber\\
&&+12(c-A)^{2}\gamma_{1}\gamma_{2}\varepsilon^{3}u^{3}-8(c-A)^{3}\gamma_{1}^{3}\varepsilon^{3}u^{3})(\frac{3c^{2}-3Ac+A^{2}}{(2c-A)(c-A)}uu_{\xi}+\frac{(c-A)B_{1}}{2c-A}\varepsilon u^{2}u_{\xi}\nonumber\\
&&+\frac{(c-A)B_{2}}{2c-A}\varepsilon^{2}u^{3}u_{\xi}+\frac{(c-A)B_{3}}{2c-A}\varepsilon^{3}u^{4}u_{\xi}+\frac{(c-A)^{2}}{3(2c-A)}\mu u_{\xi\xi\xi})+O(\varepsilon^{6},\varepsilon^{3}\mu,\mu^{2})\nonumber
\end{eqnarray*}
\begin{eqnarray}\label{ufang utao}
&=&-\varepsilon^{2}u^{2}(\frac{3c^{2}-3Ac+A^{2}}{(2c-A)(c-A)}uu_{\xi}+\frac{(c-A)^{2}}{3(2c-A)}\mu u_{\xi\xi\xi}+(\frac{(c-A)B_{1}}{2c-A}-\frac{2\gamma_{1}(3c^{2}-3Ac+A^{2})}{2c-A})\nonumber\\
&&\times\varepsilon u^{2}u_{\xi}+(\frac{(c-A)B_{2}}{2c-A}-\frac{2\gamma_{1}B_{1}(c-A)^{2}}{2c-A}-\frac{3\gamma_{2}(3c^{2}-3Ac+A^{2})}{2c-A}+\frac{4(c-A)\gamma_{1}^{2}}{2c-A}(3c^{2}\nonumber\\
&&-3Ac+A^{2}))\varepsilon^{2}u^{3}u_{\xi}+(\frac{(c-A)B_{3}}{2c-A}-\frac{2(c-A)^{2}B_{2}\gamma_{1}}{2c-A}-\frac{3(c-A)^{2}B_{1}\gamma_{2}}{2c-A}\nonumber\\
&&+\frac{4(c-A)^{3}\gamma_{1}^{2}B_{1}}{2c-A}-\frac{4\gamma_{3}(3c^2-3Ac+A^2)}{2c-A}+\frac{12(c-A)\gamma_{1}\gamma_{2}(3c^{2}-3Ac+A^{2})}{2c-A}\nonumber\\
&&-\frac{8(c-A)^{2}\gamma_{1}^{3}(3c^{2}-3Ac+A^{2})}{2c-A})\varepsilon^{3}u^{4}u_{\xi})+O(\varepsilon^{6},\varepsilon^{3}\mu,\mu^{2}).
\end{eqnarray}}
Hence, we get
\begin{align*}
&\varepsilon^{3}u^{3}u_{\tau}=-\varepsilon^{3}u^{3}(\frac{3c^{2}-3Ac+A^{2}}{(2c-A)(c-A)}uu_{\xi}+(\frac{(c-A)B_{1}-2\gamma_{1}(3c^{2}-3Ac+A^{2})}{2c-A})\varepsilon u^{2}u_{\xi}+(\frac{(c-A)B_{2}}{2c-A}\\
&-\frac{2\gamma_{1}B_{1}(c-A)^{2}+3\gamma_{2}(3c^{2}-3Ac+A^{2})}{2c-A}+\frac{4(c-A)\gamma_{1}^{2}(3c^{2}-3Ac+A^{2})}{2c-A})\varepsilon^{2}u^{3}u_{\xi})+O(\varepsilon^{6},\varepsilon^{3}\mu,\mu^{2}),\\
&\varepsilon^{4}u^{4}u_{\tau}=-\varepsilon^{4}u^{4}(\frac{3c^{2}-3Ac+A^{2}}{(2c-A)(c-A)}uu_{\xi}+(\frac{(c-A)B_{1}}{2c-A}-\frac{2\gamma_{1}(3c^{2}-3Ac+A^{2})}{2c-A})\varepsilon u^{2}u_{\xi})\\
&+O(\varepsilon^{6},\varepsilon^{3}\mu,\mu^{2}),\\
&\varepsilon^{5}u^{5}u_{\tau}=-\varepsilon^{5}u^{5}(\frac{3c^{2}-3Ac+A^{2}}{(2c-A)(c-A)}uu_{\xi})+O(\varepsilon^{6},\varepsilon^{3}\mu,\mu^{2}),
\end{align*}
\begin{eqnarray}\label{ugenggaocifang utao}
\varepsilon\mu u_{\tau\xi\xi}&=&-\varepsilon\mu(\frac{3c^{2}-3Ac+A^{2}}{(2c-A)(c-A)}(uu_{\xi})_{\xi\xi}+(\frac{(c-A)B_{1}}{2c-A}-\frac{2\gamma_{1}(3c^{2}-3Ac+A^{2})}{2c-A})\varepsilon (u^{2}u_{\xi})_{\xi\xi})\nonumber\\
&&+O(\varepsilon^{6},\varepsilon^{3}\mu,\mu^{2}),\nonumber\\
\varepsilon^{2}\mu u_{\xi}u_{\tau\xi}&=&-\varepsilon^{2}\mu u_{\xi}(\frac{3c^{2}-3Ac+A^{2}}{(2c-A)(c-A)}(uu_{\xi})_{\xi})+O(\varepsilon^{6},\varepsilon^{3}\mu,\mu^{2}),\nonumber\\
\varepsilon^{2}\mu u_{\tau}u_{\xi\xi}&=&-\varepsilon^{2}\mu u_{\xi\xi}(\frac{3c^{2}-3Ac+A^{2}}{(2c-A)(c-A)}uu_{\xi})+O(\varepsilon^{6},\varepsilon^{3}\mu,\mu^{2}),\nonumber\\
\varepsilon^{2}\mu uu_{\tau\xi\xi}&=&-\varepsilon^{2}\mu u(\frac{3c^{2}-3Ac+A^{2}}{(2c-A)(c-A)}(uu_{\xi})_{\xi\xi})+O(\varepsilon^{6},\varepsilon^{3}\mu,\mu^{2}).
\end{eqnarray}
Taking account of (\ref{ufang utao})-(\ref{ugenggaocifang utao}), it remains us to compute $\varepsilon uu_{\tau}$. To this end, we infer from (\ref{guanyu u}) that
{\setlength\arraycolsep{2pt}
\begin{eqnarray*}
\varepsilon uu_{\tau}&=&-\varepsilon u(2(c-A)\gamma_{1}\varepsilon uu_{\tau}+3(c-A)\gamma_{2}\varepsilon^{2}u^{2}u_{\tau}+4(c-A)\gamma_{3}\varepsilon^{3}u^{3}u_{\tau}+5(c-A)\gamma_{4}\varepsilon^{4}u^{4}u_{\tau}\\
&&-(c-A)\gamma_{6}\varepsilon\mu\frac{3c^{2}-3Ac+A^{2}}{(2c-A)(c-A)}(uu_{\xi})_{\xi\xi}+\frac{3c^{2}-3Ac+A^{2}}{(2c-A)(c-A)}uu_{\xi}+\frac{(c-A)B_{1}}{2c-A}\varepsilon u^{2}u_{\xi}\\
&&+\frac{(c-A)B_{2}}{2c-A}\varepsilon^{2}u^{3}u_{\xi}+\frac{(c-A)B_{3}}{2c-A}\varepsilon^{3}u^{4}u_{\xi}+\frac{(c-A)B_{4}}{2c-A}\varepsilon^{4}u^{5}u_{\xi}+\frac{(c-A)B_{6}}{2c-A}\varepsilon\mu u_{\xi}u_{\xi\xi}\\
&&+\frac{(c-A)B_{7}}{2c-A}\varepsilon\mu uu_{\xi\xi\xi}+\frac{(c-A)^{2}}{3(2c-A)}\mu u_{\xi\xi\xi})+O(\varepsilon^{6},\varepsilon^{3}\mu,\mu^{2}),
\end{eqnarray*}
which implies
{\setlength\arraycolsep{2pt}
\begin{align*}
&\varepsilon u(1+2(c-A)\gamma_{1}\varepsilon u+3(c-A)\gamma_{2}\varepsilon^{2}u^{2}+4(c-A)\gamma_{3}\varepsilon^{3}u^{3}+5(c-A)\gamma_{4}\varepsilon^{4}u^{4})u_{\tau}\\
=&-\varepsilon u(\frac{3c^{2}-3Ac+A^{2}}{(2c-A)(c-A)}uu_{\xi}+\frac{(c-A)B_{1}}{2c-A}\varepsilon u^{2}u_{\xi}+\frac{(c-A)B_{2}}{2c-A}\varepsilon^{2}u^{3}u_{\xi}+\frac{(c-A)B_{3}}{2c-A}\varepsilon^{3}u^{4}u_{\xi}\\
&+\frac{(c-A)B_{4}}{2c-A}\varepsilon^{4}u^{5}u_{\xi}+\frac{(c-A)^{2}}{3(2c-A)}\mu u_{\xi\xi\xi}+\frac{(c-A)B_{6}}{2c-A}\varepsilon\mu u_{\xi}u_{\xi\xi}+\frac{(c-A)B_{7}}{2c-A}\varepsilon\mu uu_{\xi\xi\xi}\\
&-\frac{\gamma_{6}(3c^{2}-3Ac+A^{2})}{2c-A}\varepsilon\mu(uu_{\xi})_{\xi\xi})+O(\varepsilon^{6},\varepsilon^{3}\mu,\mu^{2}),
\end{align*}
It thus follows that
{\setlength\arraycolsep{2pt}
\begin{eqnarray*}
&&\varepsilon uu_{\tau}=-\varepsilon u(1-(2(c-A)\gamma_{1}\varepsilon u+3(c-A)\gamma_{2}\varepsilon^{2}u^{2}+4(c-A)\gamma_{3}\varepsilon^{3}u^{3}+5(c-A)\gamma_{4}\varepsilon^{4}u^{4})+4(c\\
&& -A)^{2} \gamma_{1}^{2}\varepsilon^{2}u^{2}+12(c-A)^{2}\gamma_{1}\gamma_{2}\varepsilon^{3}u^{3}+9(c-A)^{2}\gamma_{2}^{2}\varepsilon^{4}u^{4}+16(c-A)^{2}\gamma_{1}\gamma_{3}\varepsilon^{4}u^{4}-(8(c-A)^{3}\\
&&\times \gamma_{1}^{3}\varepsilon^{3}u^{3}+36(c-A)^{3}\gamma_{1}^{2}\gamma_{2}\varepsilon^{4}u^{4})+16(c-A)^{4}\gamma_{1}^{4}\varepsilon^{4}u^{4})(\frac{3c^{2}-3Ac+A^{2}}{(2c-A)(c-A)}uu_{\xi}+\frac{(c-A)B_{1}}{2c-A}\varepsilon u^{2}u_{\xi}\\
&&+\frac{(c-A)B_{2}}{2c-A}\varepsilon^{2}u^{3}u_{\xi}+\frac{(c-A)B_{3}}{2c-A}\varepsilon^{3}u^{4}u_{\xi}+\frac{(c-A)B_{4}}{2c-A}\varepsilon^{4}u^{5}u_{\xi}+\frac{(c-A)^{2}}{3(2c-A)}\mu u_{\xi\xi\xi}+\frac{(c-A)B_{6}}{2c-A}\\
&&\times\varepsilon\mu u_{\xi}u_{\xi\xi}+\frac{(c-A)B_{7}}{2c-A}\varepsilon\mu uu_{\xi\xi\xi}-\frac{\gamma_{6}(3c^{2}-3Ac+A^{2})}{2c-A}\varepsilon\mu(uu_{\xi})_{\xi\xi})+O(\varepsilon^{6},\varepsilon^{3}\mu,\mu^{2}),
\end{eqnarray*}
and then
{\setlength\arraycolsep{2pt}
\begin{eqnarray*}
\varepsilon uu_{\tau}&=&-\varepsilon u(\frac{3c^{2}-3Ac+A^{2}}{(2c-A)(c-A)}uu_{\xi}+(\frac{(c-A)B_{1}}{2c-A}-\frac{2\gamma_{1}(3c^{2}-3Ac+A^{2})}{2c-A})\varepsilon u^{2}u_{\xi}\nonumber\\
&&+(\frac{(c-A)B_2}{2c-A}-\frac{3\gamma_2(3c^2-3Ac+A^2)}{2c-A}-\frac{2(c-A)^2\gamma_1B_1}{2c-A}\nonumber\\
&&+\frac{4(c-A)\gamma_1^2(3c^2-3Ac+A^2)}{2c-A})\varepsilon^2u^3u_{\xi}+(\frac{(c-A)B_{3}}{2c-A}-\frac{2(c-A)^{2}\gamma_{1}B_{2}}{2c-A}\nonumber\\
&&-\frac{3(c-A)^{2}B_{1}\gamma_{2}}{2c-A}-\frac{4\gamma_{3}(3c^{2}-3Ac+A^{2})}{2c-A}+\frac{4(c-A)^{3}B_{1}\gamma_{1}^{2}}{2c-A}\nonumber\\
&&+\frac{12(c-A)\gamma_{1}\gamma_{2}(3c^{2}-3Ac+A^{2})}{2c-A}-\frac{8(c-A)^{2}\gamma_{1}^{3}(3c^{2}-3Ac+A^{2})}{2c-A})\varepsilon^{3}u^{4}u_{\xi}\nonumber\\
&&+(\frac{(c-A)B_{4}}{2c-A}-\frac{2(c-A)^{2}\gamma_{1}B_{3}}{2c-A}-\frac{3(c-A)^{2}\gamma_{2}B_{2}}{2c-A}-\frac{4(c-A)^{2}\gamma_{3}B_{1}}{2c-A}\nonumber\\
&&-\frac{5\gamma_{4}(3c^{2}-3Ac+A^{2})}{2c-A}+\frac{4(c-A)^{3}\gamma_{1}^{2}B_{2}}{2c-A}+\frac{12(c-A)^{3}\gamma_{1}\gamma_{2}B_{1}}{2c-A}\nonumber\\
&&+\frac{9(c-A)\gamma_{2}^{2}(3c^{2}-3Ac+A^{2})}{2c-A}+\frac{16(c-A)\gamma_{1}\gamma_{3}(3c^{2}-3Ac+A^{2})}{2c-A}-\frac{8(c-A)^{4}\gamma_{1}^{3}B_{1}}{2c-A}\nonumber\\
&&-\frac{36(c-A)^{2}\gamma_{1}^{2}\gamma_{2}(3c^{2}-3Ac+A^{2})}{2c-A}+\frac{16(c-A)^{3}\gamma_{1}^{4}(3c^{2}-3Ac+A^{2})}{2c-A})\varepsilon^{4}u^{5}u_{\xi}\nonumber
\end{eqnarray*}
\begin{eqnarray}\label{u utao}
&&+\frac{(c-A)^{2}}{3(2c-A)}\mu u_{\xi\xi\xi}+\frac{(c-A)B_{6}}{2c-A}\varepsilon\mu u_{\xi}u_{\xi\xi}+(\frac{(c-A)B_{7}}{2c-A}-\frac{2(c-A)^3\gamma_1}{3(2c-A)})\varepsilon\mu uu_{\xi\xi\xi}\nonumber\\
&&-\frac{\gamma_{6}(3c^{2}-3Ac+A^{2})}{2c-A}\varepsilon\mu(uu_{\xi})_{\xi\xi})+O(\varepsilon^{6},\varepsilon^{3}\mu,\mu^{2}).
\end{eqnarray}
Decomposing the term $\varepsilon\mu u_{\tau\xi\xi}$ into $\varepsilon\mu(1-\nu)u_{\tau\xi\xi}+\varepsilon\mu\nu u_{\tau\xi\xi}$ for some constant $\nu$ (to be determined later),
and substituting (\ref{ufang utao})-(\ref{u utao}) into (\ref{guanyu u}), we obtain
{\setlength\arraycolsep{2pt}
\begin{eqnarray}\label{u kexi tao equation}
u_{\tau}+&&\frac{3c^{2}-3Ac+A^{2}}{(2c-A)(c-A)}uu_{\xi}+\frac{(c-A)^{2}}{3(2c-A)}\mu u_{\xi\xi\xi}+(c-A)\gamma_{6}(1-\nu)\varepsilon\mu u_{\tau\xi\xi}+B_{11}\varepsilon u^{2}u_{\xi}\nonumber\\
&&+B_{12}\varepsilon^{2}u^{3}u_{\xi}+B_{13}\varepsilon^{3}u^{4}u_{\xi}+B_{14}\varepsilon^{4}u^{5}u_{\xi}+B_{15}\varepsilon^{5}u^{6}u_{\xi}+\varepsilon\mu(B_{16}u_{\xi}u_{\xi\xi}+B_{17}uu_{\xi\xi\xi})\nonumber\\
&&+\varepsilon^{2}\mu(B_{18}u^{2}u_{\xi\xi\xi}+B_{19}u_{\xi}^{3}+B_{20}uu_{\xi}u_{\xi\xi})=O(\varepsilon^{6},\varepsilon^{3}\mu,\mu^{2}),
\end{eqnarray}
where $B_{11}=\frac{c^5(c^2-1)(c^2+2)}{2(c^2+1)^3},
B_{12}=\frac{c^6(c-1)^2(c+1)^2(c^4+4c^2+6)}{6(c^2+1)^5},
B_{13}=\frac{c^7(c-1)^3(c+1)^3(c^2+4)(c^4+2c^2+6)}{24(c^2+1)^7},
B_{14}\\=\frac{c^8(c-1)^4(c+1)^4(c^8+8c^6+28c^4+36c^2+120)}{120(c^2+1)^9},\
B_{15}=\frac{c^9(c-1)^5(c+1)^5(c^{10}+10c^8+44c^6+152c^4-108c^2+720)}{720(c^2+1)^{11}},\
B_{16}=\\
\frac{(c^2+c+1)(c^2-c+1)\gamma_6(1-3\nu)}{c(c^2+1)}+\frac{2c^6+7c^4+14c^2+6}{3(c^2+1)^3},
B_{17}=\frac{(c^2+c+1)(c^2-c+1)\gamma_6}{c(c^2+1)}(1-\nu)+\frac{c^4+6c^2+3}{3(c^2+1)^3},
B_{18}=\frac{\gamma_6(1-\nu)}{2(c^2+1)^5}\\
\times(c^{12}+3c^{10}+c^8-3c^6-2c^4)+\frac{(c^{10}+4c^8+9c^6+37c^4+24c^2+5)c}{6(c^2+1)^5},
B_{19}=-\frac{\gamma_6\nu}{(c^2+1)^5}(c^{12}+3c^{10}+c^8-3c^6-2c^4)-\frac{(c^{12}+5c^{10}+11c^8+14c^6+11c^4+5c^2+1)\gamma_7}{c(c^2+1)^5}
+\frac{6c^{12}+27c^{10}+53c^8+70c^6+41c^4+9c^2}{6c(c^2+1)^5}),
B_{20}=\frac{-(3c^{12}+9c^{10}+3c^8-9c^6-6c^4)\gamma_6\nu}{(c^2+1)^5}-\frac{(2c^{12}+10c^{10}+11c^8+28c^6+11c^4+10c^2+2)\gamma_8}{c(c^2+1)^5}+\frac{(3c^{12}+11c^{10}+15c^{8}+11c^6+6c^4+2c^2)\gamma_6}{(c^2+1)^5}+\frac{5c^{12}+24c^{10}+49c^8+90c^6+50c^4+11c^2}{3c(c^2+1)^5}.$

Back to the original transformation $x=\varepsilon^{-\frac{1}{2}}\xi+c\varepsilon^{-\frac{3}{2}}\tau,\ t=\varepsilon^{-\frac{3}{2}}\tau$, we have
\begin{align}
\nonumber
\frac{\partial}{\partial\xi}=\varepsilon^{-\frac{1}{2}}\partial_x, \qquad \frac{\partial}{\partial\tau}=\varepsilon^{-\frac{3}{2}}(c\partial_x+\partial_t).
\end{align}
Making use of this transformation, the equation (\ref{u kexi tao equation}) can be written as
{\setlength\arraycolsep{2pt}
\begin{align*}
u_{t}&+cu_{x}+\frac{3c^{2}-3Ac+A^{2}}{(2c-A)(c-A)}\varepsilon uu_{x}+B_{11}\varepsilon^{2}u^{2}u_{x}+B_{12}\varepsilon^{3}u^{3}u_{x}+B_{13}\varepsilon^{4}u^{4}u_{x}+B_{14}\varepsilon^{5}u^{5}u_{x}\\
&+(c-A)\gamma_{6}(1-\nu)\mu u_{txx}+(\frac{(c-A)^{2}}{3(2c-A)}+c(c-A)\gamma_{6}(1-\nu))\mu u_{xxx}+\varepsilon\mu(B_{16}u_{x}u_{xx}\\
&+B_{17}uu_{xxx})+\varepsilon^{2}\mu(B_{18}u^{2}u_{xxx}+B_{19}u_{x}^{3}+B_{20}uu_{x}u_{xx})=O(\varepsilon^{6},\varepsilon^{3}\mu,\mu^{2}).
\end{align*}
In order to obtain the CH-type terms from the above equation, it is required that
\begin{eqnarray*}
\frac{2(3c^{2}-3Ac+A^{2})}{3(2c-A)(c-A)}(c-A)\gamma_{6}(1-\nu)=2B_{17}=B_{16},
\end{eqnarray*}
which leads to
\begin{eqnarray*}
\frac{2(3c^{2}-3Ac+A^{2})}{3(2c-A)}\gamma_{6}=\frac{4c^6+7c^4-14c^2-9}{18(c^2+1)^3},
\end{eqnarray*}
and then
\begin{eqnarray*}
\frac{2(3c^{2}-3Ac+A^{2})}{3(2c-A)}\gamma_{6}(1-\nu)=2B_{17}=B_{16}=-\frac{c^4+6c^2+3}{3(c^2+1)^3}.
\end{eqnarray*}
Therefore, we establish the desired shallow-water Eq. (\ref{Vorticity}) with constant vorticity under the scaling $\mu \ll 1$, $\varepsilon=O(\sqrt[3]{\mu})$ in the form
{\setlength\arraycolsep{2pt}
\begin{eqnarray*}
u_{t}-\beta\mu u_{xxt}&+&cu_{x}+3\alpha\varepsilon uu_{x}-\beta_{0}\mu u_{xxx}+\omega_{1}\varepsilon^{2}u^{2}u_{x}+\omega_{2}\varepsilon^{3}u^{3}u_{x}+\omega_{3}\varepsilon^{4}u^{4}u_{x}+\omega_{4}\varepsilon^{5}u^{5}u_{x}\\
&=&\alpha\beta\varepsilon\mu(2u_{x}u_{xx}+uu_{xxx})+\varepsilon^{2}\mu(\omega_{5}u^{2}u_{xxx}+
\omega_{6}u^{3}_{x}+\omega_{7}uu_{x}u_{xx}),
\end{eqnarray*}
and get the height parameter $z$ in $\gamma_6,\gamma_7,\gamma_8$ may take the value $z_0=\frac{\sqrt{6(c^4+c^2+1)(2c^8+14c^6+23c^4-3)}}{6(c^4+c^2+1)(c^2+1)}.$

\begin{remark2}\label{Rem2.1}
If we apply the scaling $\mu \ll 1$, $\varepsilon=O(\sqrt[4]{\mu})$ to derive a shallow-water waves model with constant vorticity, it requires us to continue calculating the asymptotic expansion at least at order $O(\varepsilon^{4}\mu^{1}).$  According to the similar calculation as the order $O(\varepsilon^{3}\mu^{1})$, we have to obtain an expression of $u_{31}$ in terms of $\eta_{ij}$. However, this is actually not feasible unless we use an integral expression. More specifically, thanks to (\ref{eta21tau}), we get
{\setlength\arraycolsep{2pt}
\begin{eqnarray*}
\eta_{21,\tau}&=&2c_1(\eta_{21}\eta_{00}+\eta_{20}\eta_{01}+\eta_{11}\eta_{10})_{\xi}+3c_2(\eta_{11}\eta_{00}^2+2\eta_{10}\eta_{01}\eta_{00})_{\xi}+4c_3(\eta_{01}\eta_{00}^3)_{\xi}\nonumber\\
&&-\frac{(c-A)^2}{3(2c-A)}\eta_{20,\xi\xi\xi}+2c_6(\eta_{10,\xi}\eta_{00,\xi})_{\xi}+c_7(\eta_{00}\eta_{10,\xi\xi}+\eta_{10}\eta_{00,\xi\xi})_{\xi}-\frac{c_8}{2c-A}\eta_{00,\xi}^3\nonumber\\
&&-\frac{c_9}{2c-A}\eta_{00}^2\eta_{00,\xi\xi\xi}-\frac{c_{10}}{2c-A}\eta_{00}\eta_{00,\xi}\eta_{00,\xi\xi}.
\end{eqnarray*}}
Then we have
{\setlength\arraycolsep{2pt}
\begin{eqnarray*}
u_{31,\xi}&=&(c-A)\eta_{31,\xi}-2(c+c_{1}-\frac{A}{2})(\eta_{00}\eta_{21}+\eta_{10}\eta_{11}+\eta_{20}\eta_{01})_{\xi}+3(c+c_1-c_2-\frac{A}{2})(\eta^{2}_{00}\eta_{11}\\
&&+2\eta_{00}\eta_{10}\eta_{01})_{\xi}-4(c+c_1-c_2+c_3-\frac{A}{2})(\eta_{00}^{3}\eta_{01})_{\xi}+(\frac{c-A}{6}+\frac{(c-A)^2}{3(2c-A)}\\
&&-\frac{c-A}{2}z^2)\eta_{20,\xi\xi\xi}-(2c_6+\frac{2}{3}(c+c_1-\frac{A}{2})-2(c+c_1-\frac{A}{2})z^2)(\eta_{00,\xi}\eta_{10,\xi})_{\xi}\\
&&-(c_{7}+\frac{(c-A)^{2}}{3(2c-A)}+\frac{1}{3}(c_{1}+\frac{A}{2}))(\eta_{10}\eta_{00,\xi\xi}+\eta_{00}\eta_{10,\xi\xi})_{\xi}\\
&&+(\frac{c_8}{2c-A}+c_6+\frac{1}{3}c+\frac{1}{3}c_1-c_2-\frac{A}{6}-3(c+c_1-c_2-\frac{A}{2})z^2)(\eta_{00}\eta_{00,\xi}^2)_{\xi}\\
&&+(\frac{c_9}{2c-A}+c_7-\frac{1}{6}c_1-\frac{1}{2}c_2-\frac{A}{2}+\frac{(c-A)^2}{3(2c-A)}-\frac{3}{2}(c+c_1-c_2-\frac{A}{2})z^2)(\eta_{00}^2\eta_{00,\xi\xi})_{\xi}\\
&&-(\frac{2c_8}{2c-A}+\frac{2c_9}{2c-A}-\frac{c_{10}}{2c-A})\eta_{00}\eta_{00,\xi}\eta_{00,\xi\xi}.
\end{eqnarray*}}
In view of the above equality, $2c_{8}+2c_{9}-c_{10}=\frac{c^6+7c^4+7c^2+3}{3c^2(c^{2}+1)^{3}}$ is not identical to zero. Thus the term $\eta_{00}\eta_{00,\xi}\eta_{00,\xi\xi}$
can not be written as a $\xi$-derivative. Hence, when there are no non local terms appearing in the free surface equation, we may calculate the asymptotic expansion order up to the order $O(\varepsilon^{3}\mu^{1})$.
\end{remark2}

\section{Local well-posedness}
\newtheorem{theorem3}{Theorem}[section]
\newtheorem{lemma3}{Lemma}[section]
\newtheorem {remark3}{Remark}[section]
\newtheorem {definition3}{Definition}[section]
\newtheorem{corollary3}{Corollary}[section]
\par
In the sequel, for notational convenience, we shall consider the following initial value problem of Eq. (\ref{nonlocal}) with more general coefficients
\begin{equation}\label{equation}
\left\{\begin{array}{ll}u_t+(\alpha_1+\alpha_2 u+\alpha_3 u^2)u_x=\Lambda^{-2}f(u,u_x), \ & t>0,x\in \mathbb{R},\\
u(0,x)=u_0(x),  & x\in \mathbb{R},
\end{array}\right.
\end{equation}
where the operator $\Lambda:=(1-\partial_x^2)^{\frac{1}{2}}$, and the function $f(u,u_x):=\partial_x\big(\sum_{i=1}^{6}\beta_iu^i+\beta_7u^2_x+\beta_8uu^2_x\big)+\gamma u^3_x,$
$\alpha_i,\beta_i,\gamma \in \R$. In this section, we first prove the local well-posedness in Besov spaces $B_{p,r}^s$ with $s>\max({1+\frac{1}{p}},\frac{3}{2}),1\leq p, r\leq +\infty$. Then, for $s=\frac{3}{2}$, we show the local well-posedness in $B_{2,1}^{\frac{3}{2}}$. More precisely, we have
\begin{theorem3}\label{Th3.1}
Let $s,p,r$ satisfy the condition $s>\max({1+\frac{1}{p}},\frac{3}{2}), 1\leq p,r\leq +\infty,$
or $s=\frac{3}{2},p=2,r=1$. If $u_0\in
B^{s}_{p,r}$, then there exists a time $T>0$, such that the initial value problem (\ref{equation}) has a unique solution $u\in E^{s}_{p,r}(T)$.
Moreover, the solution depends continuously on
the initial data, i.e.,the map $u_0\rightarrow u$ is continuous from a neighborhood of $u_0\in B^{s}_{p,r}$ into
$C([0,T];B^{s'}_{p,r})\cap C^1([0,T];B^{s'-1}_{p,r})$
for all $s'<s$ when $r=+\infty$ and $s'=s$ whereas $r<+\infty$.
\end{theorem3}

\begin{remark3}\label{Rem3.1}
As is well-known, for every $s\in \R$, $B_{2,2}^s(\R)=H^s(\R)$. Thus Theorem \ref{Th3.1} holds
true in the corresponding Sobolev spaces with $s>\frac{3}{2}$, which recovers the well-posedness results in \cite{L-L-L} obtained by Kato's semigroup theory \cite{Kato}.
Moreover, in view of the continuous embedding for $s>\frac{3}{2}$:  $H^s(\R) \hookrightarrow B_{2,1}^{\frac{3}{2}}(\R)$, we improve our previous local well-posedness results in \cite{L-L-L}.
\end{remark3}

\subsection{Well-posedness in $B_{p,r}^s$ with $s>\max({1+\frac{1}{p}},\frac{3}{2}),1\leq p, r\leq +\infty$}\label{Sub3.1}
In this subsection, we prove Theorem \ref{Th3.1} for the case of $s>\max({1+\frac{1}{p}},\frac{3}{2}),1\leq p, r\leq +\infty$.
Motivated by \cite{Danchin}, we shall apply the classical Friedrichs regularization method to
construct the approximate solutions. The proof is similar to \cite{Danchin,Zhou}, but more nonlinear terms appear in (\ref{equation}).
Thus we shall only show how to get the uniformly boundedness of the approximate solutions. For simplicity, we omit the details of the other parts of the proof.

By a standard iterative process, starting from $u^{(0)}:=0$, we define by induction a sequence of smooth functions $\{u^{(n)}\}_{n\in \mathbb{N}}$ by solving the following linear transport equation
\begin{equation}\label{approximate}
\left\{\begin{array}{ll}u^{(n+1)}_t+\big(\alpha_1+\alpha_2 u^{(n)}+\alpha_3 (u^{(n)})^2\big)u^{(n+1)}_x=\Lambda^{-2}f(u^{(n)},u_x^{(n)}),  \ & t>0,x\in \mathbb{R}, \\
u^{(n+1)}\mid_{t=0}=u^{(n+1)}_{0}(x)=S_{n+1}u_0, & x\in \mathbb{R}.
\end{array}\right.
\end{equation}
Since all the data $S_{n+1}u_0 \in B^{\infty}_{p,r}$, it follows from Lemma \ref{Lem5.2} by induction that for all $n\in \mathbb{N}$, Eq. (\ref{approximate}) has a global solution belonging to
$C(\R^+;B^{\infty}_{p,r})$. Notice that $\Lambda^{-2}\in Op(S^{-2})$ and $\Lambda^{-2}\partial_x\in Op(S^{-1})$ in the sense of $(vii)$ in Proposition \ref{Pro5.1}.
According to $(iv)$ in Proposition \ref{Pro5.1}, we have
\begin{eqnarray}\label{3.3}
\|\Lambda^{-2}f(u^{(n)},u_x^{(n)})\|_{B^{s}_{p,r}}\leq\frac{C}{2}(\sum_{i=1}^{6}\|u^{(n)}\|^i_{B^{s}_{p,r}}),
\end{eqnarray}
and
\begin{eqnarray}\label{3.4}
\|\alpha_2\partial_x(u^{(n)})+\alpha_3\partial_x((u^{(n)})^2)\|_{B^{s-1}_{p,r}} \leq C(\|u^{(n)}\|_{B^{s}_{p,r}}+\|u^{(n)}\|^2_{B^{s}_{p,r}}).
\end{eqnarray}
Thanks to $(i)$ of Lemma \ref{Lem5.1} and (\ref{3.3})-(\ref{3.4}), we obtain
\begin{eqnarray}\label{3.5}
\|u^{(n+1)}(t)\|_{B^s_{p,r}}&\leq& \exp\big(C\int^t_0(\|u^{(n)}(t')\|_{B^{s}_{p,r}}+\|u^{(n)}(t')\|^2_{B^{s}_{p,r}})dt'\big)\|u_0\|_{B^s_{p,r}}\nonumber\\
&&+\frac{C}{2}\int^t_0\exp\big(C\int^t_\tau (\|u^{(n)}(t')\|_{B^{s}_{p,r}}+\|u^{(n)}(t')\|^2_{B^{s}_{p,r}})dt' \big)(\sum_{i=1}^{6}\|u^{(n)}\|^i_{B^{s}_{p,r}})d\tau.\quad\quad
\end{eqnarray}

Choosing $0<T\leq \frac{1}{3C\big(4\|u_0\|_{B^s_{p,r}}+\sum_{i=1,i\neq2}^{6}2^{2i-2}\|u_0\|^{i-1}_{B^s_{p,r}}\big)}$, by induction, we prove that
\begin{eqnarray}\label{3.6}
\|u^{n}(t)\|_{B^s_{p,r}}\leq\frac{2\|u_0\|_{B^s_{p,r}}}{1-4C\|u_0\|_{B^s_{p,r}}t},\quad \forall t\in [0,T].
\end{eqnarray}
Assume that (\ref{3.6}) is valid for $n$, by using the mean-value theorem of integrals, we get for $\xi\in (\tau,t)$
\begin{eqnarray}\label{3.7}
&&\exp\big(C\int^t_\tau (\|u^{(n)}(t')\|_{B^{s}_{p,r}}+\|u^{(n)}(t')\|^2_{B^{s}_{p,r}})dt' \big)\nonumber \\
&\leq&\exp\big(\int^t_\tau (\frac{2\|u_0\|_{B^s_{p,r}}}{1-4C\|u_0\|_{B^s_{p,r}}t'}
+\frac{4\|u_0\|^2_{B^s_{p,r}}}{(1-4C\|u_0\|_{B^s_{p,r}}t')^2})dt' \big)\nonumber\\
&\leq&\exp\big(\frac{1}{2}\ln\frac{1-4C\|u_0\|_{B^s_{p,r}}\tau}{1-4C\|u_0\|_{B^s_{p,r}}t}+
\frac{4\|u_0\|^2_{B^s_{p,r}}(t-\tau)}{(1-4C\|u_0\|_{B^s_{p,r}}\xi)^2}\big)\nonumber\\
&\leq&\exp\big(16C\|u_0\|^2_{B^s_{p,r}}T\big)\big(\frac{1-4C\|u_0\|_{B^s_{p,r}}\tau}{1-4C\|u_0\|_{B^s_{p,r}}t}\big)^{\frac{1}{2}}\leq e^{\frac{1}{3}}\big(\frac{1-4C\|u_0\|_{B^s_{p,r}}\tau}{1-4C\|u_0\|_{B^s_{p,r}}t}\big)^{\frac{1}{2}}.
\end{eqnarray}
It then follows from (\ref{3.6})-(\ref{3.7}) and the mean-value theorem of integrals again that
\begin{eqnarray}
&&\frac{C}{2}\int^t_0\exp\big(C\int^t_\tau (\|u^{(n)}(t')\|_{B^{s}_{p,r}}+\|u^{(n)}(t')\|^2_{B^{s}_{p,r}})dt' \big)(\sum_{i=1}^{6}\|u^{(n)}\|^i_{B^{s}_{p,r}})d\tau\nonumber \\
&\leq&\frac{e^{\frac{1}{3}}\|u_0\|_{B^s_{p,r}}}{(1-4C\|u_0\|_{B^s_{p,r}}t)^{\frac{1}{2}}}\int^t_0\big(\sum_{i=1}^{6}\frac{2^{i-1}C \|u_0\|^{i-1}_{B^s_{p,r}} }
{(1-4C\|u_0\|_{B^s_{p,r}}\tau)^{i-\frac{1}{2}}}\big)d\tau\nonumber \\
&\leq&\frac{e^{\frac{1}{3}}\|u_0\|_{B^s_{p,r}}}{(1-4C\|u_0\|_{B^s_{p,r}}t)^{\frac{1}{2}}}\Big(
\int^t_0\frac{2C\|u_0\|_{B^s_{p,r}}}{(1-4C\|u_0\|_{B^s_{p,r}}\tau)^{\frac{3}{2}}}d\tau
+\sum_{i=1,i\neq2}^{6}\frac{2^{i-1}C \|u_0\|^{i-1}_{B^s_{p,r}}t }{(1-4C\|u_0\|_{B^s_{p,r}}\xi_i)^{i-\frac{1}{2}}}\Big)\nonumber
\end{eqnarray}
\begin{eqnarray}\label{3.8}
&\leq&\frac{e^{\frac{1}{3}}\|u_0\|_{B^s_{p,r}}}{(1-4C\|u_0\|_{B^s_{p,r}}t)^{\frac{1}{2}}} \big(\frac{1}{(1-4C\|u_0\|_{B^s_{p,r}}t)^{\frac{1}{2}}}-1\big)\nonumber \\
&&+\frac{e^{\frac{1}{3}}\|u_0\|_{B^s_{p,r}}}{(1-4C\|u_0\|_{B^s_{p,r}}t)^{\frac{1}{2}}} \big(
\sum_{i=1,i\neq2}^{6}\frac{2^{i-1}C \|u_0\|^{i-1}_{B^s_{p,r}}T}{(1-4C\|u_0\|_{B^s_{p,r}}t)^{i-\frac{1}{2}}}\big).
\end{eqnarray}
where $0<\xi_i<t,\ (1\leq i\leq6,i\neq2)$. Inserting (\ref{3.7}) with $\tau=0$ and (\ref{3.8}) into (\ref{3.5}), it yields
\begin{eqnarray}\label{3.9}
\|u^{(n+1)}\|_{B^s_{p,r}}&\leq&\frac{e^{\frac{1}{3}}\|u_0\|_{B^s_{p,r}}}{1-4C\|u_0\|_{B^s_{p,r}}t}\big(1+
\sum_{i=1,i\neq2}^{6}\frac{2^{i-1}C \|u_0\|^{i-1}_{B^s_{p,r}}T}{(1-4C\|u_0\|_{B^s_{p,r}}t)^{i-1}}\big)\nonumber\\
&\leq&\frac{e^{\frac{1}{3}}\|u_0\|_{B^s_{p,r}}}{1-4C\|u_0\|_{B^s_{p,r}}t}\big(1+CT\sum_{i=1,i\neq2}^{6}2^{2i-2}\|u_0\|^{i-1}_{B^s_{p,r}}\big)\nonumber\\
&\leq&\frac{4e^{\frac{1}{3}}\|u_0\|_{B^s_{p,r}}}{3(1-4C\|u_0\|_{B^s_{p,r}}t)}\leq\frac{2\|u_0\|_{B^s_{p,r}}}{1-4C\|u_0\|_{B^s_{p,r}}t},
\end{eqnarray}
which proves that $\{u^{(n)}\}_{n\in \mathbb{N}}$ is uniformly bounded in
$C([0,T];{B^s_{p,r}}).$ Taking account of Eq. (\ref{approximate}), one can easily show that $\{\partial_tu^{(n)}\}_{n\in \mathbb{N}}$ is uniformly bounded in
$C([0,T];{B^{s-1}_{p,r}}).$ Therefore, this completes the proof of the uniform bound of $\{u^{(n)}\}_{n\in \mathbb{N}}$ in $E^s_{p,r}(T)$.

\subsection{Well-posedness in $B_{2,1}^{\frac{3}{2}}$}
In this subsection, we shall complete the proof of Theorem \ref{Th3.1} for $u_0\in B_{2,1}^{\frac{3}{2}}$. By using an argument similar to the one in Subsection \ref{Sub3.1}, one can prove there exists a time $T$ such that the initial value problem (\ref{equation}) has a solution $u \in E^\frac{3}{2}_{2,1}(T)$. Moreover, when $u_0\in B_{2,1}^{\frac{3}{2}}$, thanks to a similar proof as (\ref{3.6}), we can obtain
\begin{eqnarray}\label{3.0}
\|u\|_{B_{2,1}^{\frac{3}{2}}}\leq 4\|u_0\|_{B_{2,1}^{\frac{3}{2}}}.
\end{eqnarray}

Next we are devoted to give a priori estimate, which implies
the uniqueness of the solution. Suppose that $u,v\in
L^\infty(0,T;B^{\frac{3}{2}}_{2,\infty}\cap \mbox{Lip})\cap
C([0,T];B^\frac{1}{2}_{2,\infty})$ are two given solutions to
(\ref{equation}) with the initial data
$u_0,v_0 \in B^{\frac{3}{2}}_{2,\infty}\cap
\mbox{Lip},$ respectively. Denote $w:=u-v$ and
$w_0:=u_0-v_0$. Assume that there exists a positive
constant $C$ such that for some $T^*\leq T$,
\begin{eqnarray}\label{3.10}
\sup_{t\in[0,T^*]}\Big(e^{-C\int_0^tU(\tau)d\tau}\|w\|_{B^\frac{1}{2}_{2,\infty}} \Big)\leq1,
\end{eqnarray}
where $U(t):=\|\partial_x(\alpha_2 u+\alpha_3 u^2)\|_{B^\frac{1}{2}_{2,\infty}\cap
L^\infty}.$ We shall prove for all $t\in[0,T^*]$ that
\begin{eqnarray}\label{3.11}
\frac{\|w(t)\|_{B^\frac{1}{2}_{2,\infty}}}{e} \leq
e^{C\int_0^tU(\tau)d\tau}\Big(\frac{\|w_0\|_{B^\frac{1}{2}_{2,\infty}}}{e}
\Big)^{\exp(-C\int_0^tV(\tau)\ln(e+V(\tau))d\tau)},
\end{eqnarray}
where $V(t):=(1+\|u\|^5_{B^{\frac{3}{2}}_{2,\infty}\cap
\mbox{Lip}}+\|v\|^5_{B^{\frac{3}{2}}_{2,\infty}\cap
\mbox{Lip}})$. In particular, if
\begin{eqnarray*}
\|w_0\|_{B^\frac{1}{2}_{2,\infty}}\leq
e^{1-\exp(C\int_0^TV(t)\ln(e+V(t))dt)},
\end{eqnarray*}
then (\ref{3.11}) holds true on $[0,T],$ since the above inequality
implies that (\ref{3.10}) holds with $T^*=T.$

Indeed, $w$ solves the following initial value problem
\begin{equation*}
\left\{\begin{array}{ll}w_t+(\alpha_1+\alpha_2 u+\alpha_3 u^2)w_x=
-(\alpha_2+\alpha_3(u+v))wv_x+\Lambda^{-2}\partial_x(\tilde{f}_{1}(t,x))+\Lambda^{-2}(\tilde{f}_{2}(t,x))
,\\
w\big|_{t=0}=u_0-v_0,
\end{array}\right.
\end{equation*}
where $\tilde{f}_{1}(t,x):=\sum_{i=1}^{6}\beta_iw(\sum_{j=1}^{i}u^{i-j}v^{j-1})+\beta_7w_x(u_x+v_x)+\beta_8wv^2_x+\beta_8uw_x(u_x+v_x)$ and
$\tilde{f}_{2}(t,x):=\gamma w_x(u^2_x+u_xv_x+v^2_x)$. Note that $\Lambda^{-2}\in Op(S^{-2})$ and $\Lambda^{-2}\partial_x\in Op(S^{-1})$. It thus follows from $(iv)$ and $(ix)$ in Proposition \ref{Pro5.1} that
\begin{eqnarray}\label{3.12}
\|(\alpha_2+\alpha_3(u+v))wv_x\|_{B^\frac{1}{2}_{2,\infty}}\leq C
\|w\|_{B^{\frac{1}{2}}_{2,1}}\|v\|_{B^{\frac{3}{2}}_{2,\infty}\cap
\mbox{Lip}}\big(1+\|u\|_{B^{\frac{3}{2}}_{2,\infty}\cap
\mbox{Lip}}+\|v\|_{B^{\frac{3}{2}}_{2,\infty}\cap
\mbox{Lip}}\big),\quad
\end{eqnarray}
\begin{eqnarray}\label{3.13}
\|\Lambda^{-2}\partial_x(\tilde{f}_{1}(t,x)\|_{B^\frac{1}{2}_{2,\infty}}\leq C
\|w\|_{B^{\frac{1}{2}}_{2,1}}\big(\sum_{i=1}^{6}(\sum_{j=1}^{i}\|u\|_{B^{\frac{3}{2}}_{2,\infty}\cap
\mbox{Lip}}^{i-j}\|v\|_{B^{\frac{3}{2}}_{2,\infty}\cap
\mbox{Lip}}^{j-1})\big),
\end{eqnarray}
and
\begin{eqnarray}\label{3.14}
\|\Lambda^{-2}(\tilde{f}_{2}(t,x)\|_{B^\frac{1}{2}_{2,\infty}}\leq C
\|w\|_{B^{\frac{1}{2}}_{2,1}}\big(\|u\|_{B^{\frac{3}{2}}_{2,\infty}\cap
\mbox{Lip}}^{2}
+\|u\|_{B^{\frac{3}{2}}_{2,\infty}\cap
\mbox{Lip}}\|v\|_{B^{\frac{3}{2}}_{2,\infty}\cap
\mbox{Lip}}+\|v\|_{B^{\frac{3}{2}}_{2,\infty}\cap
\mbox{Lip}}^{2}\big).\quad\quad
\end{eqnarray}
Thus, by (\ref{3.12})-(\ref{3.14}) and Lemma \ref{Lem5.1} $(i)$, we obtain
\begin{eqnarray}\label{3.15}
&&e^{-C\int_0^tU(\tau)d\tau}\|w(t)\|_{B^\frac{1}{2}_{2,\infty}}\nonumber\\ &\leq&\|w_0\|_{B^\frac{1}{2}_{2,\infty}}
+\int^t_0 e^{-C\int_0^\tau U(t')dt'}\|-(\alpha_2+\alpha_3(u+v))wv_x+\Lambda^{-2}(\partial_x(\tilde{f}_{1}(t,x))
+(\tilde{f}_{2}(t,x)))\|_{B^\frac{1}{2}_{2,\infty}}d\tau\nonumber\\
&\leq&\|w_0\|_{B^\frac{1}{2}_{2,\infty}}+C\int^t_0 e^{-C\int_0^\tau U(t')dt'}\|w\|_{B^{\frac{1}{2}}_{2,1}}
\big(\sum_{i=1}^{6}(\sum_{j=1}^{i}\|u\|_{B^{\frac{3}{2}}_{2,\infty}\cap
\mbox{Lip}}^{i-j}\|v\|_{B^{\frac{3}{2}}_{2,\infty}\cap
\mbox{Lip}}^{j-1})\big)d\tau\nonumber\\
&\leq&\|w_0\|_{B^\frac{1}{2}_{2,\infty}}+
C\int^t_0 e^{-C\int_0^\tau U(t')dt'}
\|w\|_{B^\frac{1}{2}_{2,\infty}}\ln
\big(e+\frac{\|w\|_{B^\frac{3}{2}_{2,\infty}}}{\|w\|_{B^\frac{1}{2}_{2,\infty}}}\big)\cdot V(\tau)d\tau\nonumber\\
&\leq&\|w_0\|_{B^\frac{1}{2}_{2,\infty}}+
C\int^t_0 e^{-C\int_0^\tau U(t')dt'}
\|w\|_{B^\frac{1}{2}_{2,\infty}}\ln \big(e+\frac{V(\tau)}
{e^{-C\int_0^\tau U(t')dt'}\|w\|_{B^\frac{1}{2}_{2,\infty}}}\big)\cdot V(\tau)d\tau,
\end{eqnarray}
where we used Proposition \ref{Pro5.1} $(x)$ and Young inequality in the third inequality. Denote $W(t):=
e^{-C\int_0^tU(\tau)d\tau}\|w(t)\|_{B^\frac{1}{2}_{2,\infty}}$. Hence we can rewrite (\ref{3.15}) as follows
\begin{eqnarray*}
W(t)\leq W(0)+C\int^t_0W(\tau)\ln(e+\frac{V(\tau)}{W(\tau)})\cdot V(\tau)d\tau.
\end{eqnarray*}
Applying the fact that $\ln
(e+\frac{\alpha}{x})\leq \ln (e+\alpha)(1-\ln x)$ for $\forall \
x\in(0,1], \ \alpha>0,$ to the above inequality, and in view of (\ref{3.10}), we have
\begin{eqnarray}\label{3.16}
W(t)\leq W(0)+C\int^t_0V(\tau)\ln(e+V(\tau))\cdot W(\tau)(1-\ln W(\tau) )d\tau.
\end{eqnarray}
An application of Lemma \ref{Lem5.4} with $\mu(r)=r(1-\ln r)$ to (\ref{3.16}) under the hypothesis (\ref{3.10}) implies
\begin{eqnarray*}
\frac{W(t)}{e}\leq
\big(\frac{W(0)}{e}\big)^{\exp(-C\int^t_0V(\tau)\ln(e+V(\tau))d\tau)},
\end{eqnarray*}
which is the desired inequality (\ref{3.11}).

Now it suffices to prove that the solution depends continuously on the initial data $u_0\in B^\frac{3}{2}_{2,1}$, $i.e.$ the map $u_0\rightarrow u(\cdot,u_0)$: $B^\frac{3}{2}_{2,1}\rightarrow C([0,T];B^\frac{3}{2}_{2,1})$ is continuous. We first show the continuity in $C([0,T];B^\frac{1}{2}_{2,1})$. For a fixed $u_0\in B^\frac{3}{2}_{2,1}$ and $\delta>0,$ we claim that there exists a $T>0$ and $L>0$ such that for any $\tilde{u}_0\in B^\frac{3}{2}_{2,1}$ with $\|u_0-\tilde{u}_0\|\leq \delta$, the solution $\tilde{u}(\cdot,\tilde{u}_0)\in  C([0,T];B^\frac{3}{2}_{2,1})$ of (\ref{equation}) associated to $\tilde{u}_0$ satisfies
\begin{eqnarray}\label{3.18}
\|\tilde{u}\|_{L^\infty(0,T;B^\frac{3}{2}_{2,1})}\leq L.
\end{eqnarray}
Indeed, it follows from a similar argument as (\ref{3.6}) that
\begin{eqnarray*}
\|\tilde{u}\|_{B^\frac{3}{2}_{2,1}}\leq\frac{2\|\tilde{u}_0\|_{B^\frac{3}{2}_{2,1}}}{1-4C\|\tilde{u}_0\|_{B^\frac{3}{2}_{2,1}}t},\quad \forall t\in [0,T].
\end{eqnarray*}
where $0<T\leq \frac{1}{3C\big(4\|\tilde{u}_0\|_{B^\frac{3}{2}_{2,1}}+\sum_{i=1,i\neq2}^{6}2^{2i-2}\|\tilde{u}_0\|^{i-1}_{B^\frac{3}{2}_{2,1}}\big)}$. Thus one can choose
\begin{eqnarray*}
T=\frac{1}{8C\big(4(\|u_0\|_{B^\frac{3}{2}_{2,1}}+\delta)+\sum_{i=1,i\neq2}^{6}2^{2i-2}(\|u_0\|_{B^\frac{3}{2}_{2,1}}+\delta)^{i-1}\big)} \quad
\mbox{and} \quad L=4(\|u_0\|_{B^\frac{3}{2}_{2,1}}+\delta).
\end{eqnarray*}
Combining the a priori estimate (\ref{3.11}), and the uniform bounds (\ref{3.0}) and (\ref{3.18}), we obtain
\begin{eqnarray}\label{3.19}
\frac{\|\tilde{u}-u\|_{L^\infty(0,T;B^\frac{1}{2}_{2,\infty})}}{e} \leq
e^{C(L+L^2)T}\Big(\frac{\|\tilde{u}_0-u_0\|_{B^\frac{1}{2}_{2,\infty}}}{e}
\Big)^{\exp(-C(1+2L^5)\ln(e+(1+2L^5))T)},
\end{eqnarray}
provided that
$\|\tilde{u}_0-u_0\|_{B^\frac{1}{2}_{2,\infty}}\leq e^{1-\exp(C(1+2L^5)\ln(e+(1+2L^5))T)}.$ In view of the interpolation inequality $(vi)$ in Proposition \ref{Pro5.1}, the uniform bounds (\ref{3.0}) and (\ref{3.18}), and (\ref{3.19}), we deduce that the map $u_0\rightarrow u(\cdot,u_0)$: $B^\frac{3}{2}_{2,1}\rightarrow C([0,T];B^\frac{1}{2}_{2,1})$ is continuous.

Then we show the continuity in $C([0,T];B^\frac{3}{2}_{2,1})$. Assume that
$u^{(\infty)}_0\in B^\frac{3}{2}_{2,1}$ and $\{u_0^{(n)}\}_{n\in
\mathbb{N}}$ tend to $u^{(\infty)}_0$ in $B^\frac{3}{2}_{2,1}$.
Let $u^{(n)}$ be the solution of (\ref{equation})
corresponding to datum $u_0^{(n)}.$ According to the above result of the continuity in $C([0,T];B^\frac{1}{2}_{2,1})$,
proving that $u^{(n)}$ tends to $u^{(\infty)}$ in $C([0,T];B^\frac{3}{2}_{2,1})$ amounts to proving that
$v^{(n)}:= \partial_xu^{(n)}$ tends to $v^{(\infty)}:= \partial_xu^{(\infty)}$ in $C([0,T];B^\frac{1}{2}_{2,1})$.
We find that $v^{(n)}$ solves the following linear transport equation
\begin{equation*}
\left\{\begin{array}{ll}\partial_tv^{(n)}+\big(\alpha_1+\alpha_2 u^{(n)}+\alpha_3 (u^{(n)})^2\big)\partial_xv^{(n)}=\tilde{f}^{(n)}(t,x)
,\  & t>0,x\in \mathbb{R},\\
v^{(n)}\big|_{t=0}= \partial_xu_0^{(n)}(x), \ & x\in \mathbb{R},
\end{array}\right.
\end{equation*}
where $\tilde{f}^{(n)}(t,x):=-\big(\alpha_2(u^{(n)}_x)^2+2\alpha_3u^{(n)}(u_x^{(n)})^2\big)
+\Lambda^{-2}\big(\sum_{i=1}^{6}\beta_i(u^{(n)})^i+\beta_7(u^{(n)}_x)^2+\beta_8u^{(n)}(u^{(n)}_x)^2\big)\\
-\big(\sum_{i=1}^{6}\beta_i(u^{(n)})^i+\beta_7(u^{(n)}_x)^2+\beta_8u^{(n)}(u^{(n)}_x)^2\big)+\Lambda^{-2}\partial_x
\big(\gamma (u^{(n)}_x)^3\big).$ Following the method in \cite{Kato}, we decompose $v^{(n)}$ into the sum of $v_1^{(n)}$ and $v_2^{(n)}$ with
\begin{equation}\label{3.20}
\left\{\begin{array}{ll}\partial_tv_1^{(n)}+\big(\alpha_1+\alpha_2 u^{(n)}+\alpha_3 (u^{(n)})^2\big)\partial_xv_1^{(n)}=\tilde{f}^{(n)}(t,x)-\tilde{f}^{(\infty)}(t,x)
,\  & t>0,x\in \mathbb{R},\\
v_1^{(n)}\big|_{t=0}= \partial_xu_0^{(n)}(x)-\partial_xu_0^{(\infty)}(x), \ & x\in \mathbb{R}.
\end{array}\right.
\end{equation}
and
\begin{equation}\label{3.21}
\left\{\begin{array}{ll}\partial_tv_2^{(n)}+\big(\alpha_1+\alpha_2 u^{(n)}+\alpha_3 (u^{(n)})^2\big)\partial_xv_2^{(n)}=\tilde{f}^{(\infty)}(t,x)
,\  & t>0,x\in \mathbb{R},\\
v_2^{(n)}\big|_{t=0}= \partial_xu_0^{(\infty)}(x), \ & x\in \mathbb{R},
\end{array}\right.
\end{equation}
According to $\Lambda^{-2}\in Op(S^{-2})$, $\Lambda^{-2}\partial_x\in Op(S^{-1})$ and $(iv)$ of Proposition \ref{Pro5.1},
we have $\{\tilde{f}^{(n)}\}_{n\in
\bar{\mathbb{N}}}$ is uniformly bounded in
$C([0,T];B^\frac{1}{2}_{2,1})$. Moreover,
\begin{eqnarray*}
\|\tilde{f}^{(n)}-\tilde{f}^{(\infty)}\|_{B^\frac{1}{2}_{2,1}}
\leq C\big(1+\|u^{(n)}\|^5_{B^\frac{3}{2}_{2,1}}+\|u^{(\infty)}\|^5_{B^\frac{3}{2}_{2,1}}\big)
\big(\|u^{(n)}-u^{(\infty)}\|_{B^\frac{1}{2}_{2,1}}+\|u_x^{(n)}-u_x^{(\infty)}\|_{B^\frac{1}{2}_{2,1}}\big).
\end{eqnarray*}
Then applying Lemma \ref{Lem5.1} to (\ref{3.20}), it yields
\begin{eqnarray}\label{3.22}
\|v_1^{(n)}(t)\|_{B^\frac{1}{2}_{2,1}}&\leq&
\exp\big(C\int^t_0(\|u^{(n)}(t')\|_{B^\frac{3}{2}_{2,1}}+\|u^{(n)}(t')\|^2_{B^\frac{3}{2}_{2,1}})dt'\big)
\|\partial_xu_0^{(n)}-\partial_xu_0^{(\infty)}\|_{B^\frac{1}{2}_{2,1}} \nonumber\\
&&+C\int^t_0\exp\big(C\int^t_\tau (\|u^{(n)}(t')\|_{B^\frac{3}{2}_{2,1}}+\|u^{(n)}(t')\|^2_{B^\frac{3}{2}_{2,1}})dt' \big)
\big(1+\|u^{(n)}\|^5_{B^\frac{3}{2}_{2,1}}\nonumber\\
&&+\|u^{(\infty)}\|^5_{B^\frac{3}{2}_{2,1}}\big)\cdot
\big(\|u^{(n)}-u^{(\infty)}\|_{B^\frac{1}{2}_{2,1}}+\|u_x^{(n)}-u_x^{(\infty)}\|_{B^\frac{1}{2}_{2,1}}\big)
d\tau\nonumber\\
&\leq& e^{C(L+L^2)T}\big(\|\partial_xu_0^{(n)}-\partial_xu_0^{(\infty)}\|_{B^\frac{1}{2}_{2,1}}+C\int^t_0
(1+2L^5)\big(\|u^{(n)}-u^{(\infty)}\|_{B^\frac{1}{2}_{2,1}}\nonumber\\
&&+\|u_x^{(n)}-u_x^{(\infty)}\|_{B^\frac{1}{2}_{2,1}}\big)
d\tau\big).
\end{eqnarray}
On the other hand, note that $\{u^{(n)}\}_{n\in
\bar{\mathbb{N}}}$ is
uniformly bounded in $C([0,T];B^\frac{3}{2}_{2,1})$ and
$u^{(n)} \rightarrow u^{(\infty)}\ \mbox{in}\
C([0,T];B^\frac{1}{2}_{2,1})$ as $n\rightarrow\infty.$ Then applying Lemma \ref{Lem5.3} to (\ref{3.21}), we find that
$v_2^{(n)}\rightarrow v^{(\infty)}= \partial_xu^{(\infty)}$ in $C([0,T];B^\frac{1}{2}_{2,1})$ as $n\rightarrow\infty.$
Hence, combining  (\ref{3.22}) with the convergences $v_2^{(n)}\rightarrow u_x^{(\infty)}$ and  $u^{(n)}\rightarrow u^{(\infty)}$ in $C([0,T];B^\frac{1}{2}_{2,1})$, for sufficiently small $\varepsilon>0$ and large enough $n\in \mathbb{N}$, we obtain
\begin{eqnarray*}
\|\partial_xu^{(n)}-\partial_xu^{(\infty)}\|_{B^\frac{1}{2}_{2,1}}
&\leq& \|v_1^{(n)}\|_{B^\frac{1}{2}_{2,1}}+\|v_2^{(n)}-u_x^{(\infty)}\|_{B^\frac{1}{2}_{2,1}}\\
&\leq&\varepsilon + C_{L,T}
\big(\|\partial_xu_0^{(n)}-\partial_xu_0^{(\infty)}\|_{B^\frac{1}{2}_{2,1}}+\int^t_0
\big(\varepsilon
+\|\partial_xu^{(n)}-\partial_xu^{(\infty)}\|_{B^\frac{1}{2}_{2,1}}\big)d\tau\big).
\end{eqnarray*}
Therefore, applying Gronwall's lemma to the above inequality, we
then get the desired result of the continuity in $C([0,T];B^\frac{3}{2}_{2,1})$.

\section{Blow-up scenario}
\newtheorem{theorem4}{Theorem}[section]
\newtheorem{lemma4}{Lemma}[section]
\newtheorem {remark4}{Remark}[section]
\newtheorem {definition4}{Definition}[section]
\newtheorem{corollary4}{Corollary}[section]
\par
In this section, we will give a more precise blow-up criterion compared with our previous one given in \cite{L-L-L}. It reads
\begin{theorem4}\label{Th4.1}
Let $u_0\in B^\frac{3}{2}_{2,1}$. Assume $T^\star$ is the maximal existence time of the corresponding solution $u$ to Eq. (\ref{equation}) guaranteed by Theorem \ref{Th3.1}.
Then we have
\begin{eqnarray*}
T^\star<+\infty\Rightarrow \int_0^{T^\star} \|u_x\|^2_{L^\infty}d\tau=+\infty.
\end{eqnarray*}
\end{theorem4}
\begin{proof}
Applying $\Delta_q$ to Eq. (\ref{equation}), we have
\begin{eqnarray}\label{4.1}
\big(\partial_t+(\alpha_1+\alpha_2 u+\alpha_3 u^2)\partial_x\big)\Delta_q u=[\alpha_2 u+\alpha_3 u^2,\Delta_q]\partial_xu+\Delta_q\Lambda^{-2}f(u,u_x),
\end{eqnarray}
where $[\cdot,\cdot]$ denotes the commutator of the operators. In view of the commutator estimates of Lemma 2.100 in \cite{Bahouri}, we obtain
\begin{eqnarray}\label{4.2}
\|(2^{\frac{3}{2}q}\|[\alpha_2 u+\alpha_3 u^2,\Delta_q]\partial_xu\|_{L^2})_{q\geq -1}\|_{l^1}
&\leq& C(\|(\alpha_2 u+\alpha_3 u^2)_x\|_{L^\infty}\|u\|_{B^\frac{3}{2}_{2,1}}
+\|u_x\|_{L^\infty}\|u-u^2\|_{B^\frac{3}{2}_{2,1}})\nonumber \\
&\leq& C(1+\|u\|_{L^\infty})\|u_x\|_{L^\infty}\|u\|_{B^\frac{3}{2}_{2,1}}.
\end{eqnarray}
where we used 1-$D$ Moser-type estimate $(viii)$ in Proposition \ref{Pro5.1}. Note that $\Lambda^{-2}\in Op(S^{-2})$ and $\Lambda^{-2}\partial_x\in Op(S^{-1})$, we get
\begin{eqnarray}\label{4.3}
\|\Lambda^{-2}f(u,u_x)\|_{B^\frac{3}{2}_{2,1}}\leq C (\sum_{i=1}^{6}\|u\|^{i-1}_{L^\infty}+\|u_x\|_{L^\infty}+\|u\|_{L^\infty}\|u_x\|_{L^\infty})\|u\|_{B^\frac{3}{2}_{2,1}}.
\end{eqnarray}
Going along the lines of the proof of Lemma \ref{Lem5.1}, by using the estimates (\ref{4.2})-(\ref{4.3}), it yields
\begin{eqnarray*}
\|u\|_{B^\frac{3}{2}_{2,1}}\leq\|u_0\|_{B^\frac{3}{2}_{2,1}}+C\int_0^t\big(\sum_{i=1}^{6}\|u\|^{i-1}_{L^\infty}+\|u_x\|_{L^\infty}+\|u\|_{L^\infty}\|u_x\|_{L^\infty}\big)\|u\|_{B^\frac{3}{2}_{2,1}}d\tau.
\end{eqnarray*}
It thus follows from Gronwall's inequality and the embedding $H^1(\R)\hookrightarrow L^\infty(\R)$that
\begin{eqnarray}\label{4.4}
\|u\|_{B^\frac{3}{2}_{2,1}}&\leq& \|u_0\|_{B^\frac{3}{2}_{2,1}}\exp\big( C \int_0^t(\sum_{i=1}^{6}\|u\|^{i-1}_{L^\infty}+\|u_x\|_{L^\infty}+\|u\|_{L^\infty}\|u_x\|_{L^\infty})d\tau\big)\nonumber \\
&\leq& \|u_0\|_{B^\frac{3}{2}_{2,1}}\exp\big( C \int_0^t(\sum_{i=1}^{6}\|u\|^{i-1}_{H^1}+\|u_x\|_{L^\infty}+\|u\|_{H^1}\|u_x\|_{L^\infty})d\tau\big)\nonumber \\
&\leq&  \|u_0\|_{B^\frac{3}{2}_{2,1}}\exp\big( C \int_0^t(\sum_{i=1}^{6}\|u_0\|^{i-1}_{H^1}\exp(C(i-1)\int_0^\tau\|u_x\|^2_{L^\infty}dt')\nonumber \\
&&+\|u_x\|_{L^\infty}+\|u_0\|_{H^1}\exp(C\int_0^\tau\|u_x\|^2_{L^\infty}dt')\|u_x\|_{L^\infty})d\tau\big),
\end{eqnarray}
where we used $\|u\|_{H^1}\leq\|u_0\|_{H^1}\exp(C\int_0^t\|u_x\|^2_{L^\infty}d\tau),\ t\in[0,T^\star)$ given in Theorem 3.2 of \cite{L-L-L}.

Now if we assume that $\int_0^{T^\star} \|u_x\|^2_{L^\infty}d\tau\leq M$ for some $M>0,$ then we have $\int_0^{T^\star} \|u_x\|_{L^\infty}d\tau\leq M$, for finite $T^\star$.
Owing to (\ref{4.4}), we find $\|u(T^\star)\|_{B^\frac{3}{2}_{2,1}}<+\infty,$  which contradicts the fact that $T^\star$ is the maximal existence time of the solution. This completes the proof of
Theorem \ref{Th4.1}.
\end{proof}

\noindent\textbf{Acknowledgments} \
The work is supported by National Nature Science Foundation of China under Grant 12001528. The authors thank the anonymous
referee for helpful suggestions and comments.

\noindent\textbf{Conflict of interest} \
The authors declare that they have no conflict of interest.

\noindent\textbf{Data availability statements} \
The data that supports the findings of this study are available within the article.

\section{Appendix}
\newtheorem {remark5}{Remark}[section]
\newtheorem{theorem5}{Theorem}[section]
\newtheorem{definition5}{Definition}[section]
\newtheorem{proposition5}{Proposition}[section]
\newtheorem{lemma5}{Lemma}[section]
In this section, we introduce some basic theory of the
Littlewood-Paley decomposition and the transport equation theory on
Besov spaces for completeness. One can check \cite{Bahouri,Danchin,Danchin1} for details.
There exist two smooth radial functions $\chi(\xi)$ and $\varphi(\xi)$ valued in $[0,1]$, such that $\chi$ is supported in
$\mathcal{B}=\{\xi\in\R^d,|\xi|\leq\frac{4}{3}\}$ and $\varphi(\xi)$ is supported in $\mathcal{C}=
\{\xi\in\R^d,\frac{3}{4}\leq|\xi|\leq\frac{8}{3}\}$. Denote $\mathcal{F}$ by the Fourier transform and $\mathcal{F}^{-1}$ by its inverse.
For all $u\in \mathcal{S'}(\R^d)$($\mathcal{S'}(\R^d)$ denotes the tempered distribution spaces), the nonhomogeneous dyadic operators $\Delta_q$ and
the low frequency cut-off operator $S_q$ are defined as follows:
$ \Delta_q u=0$ for $q\leq-2$, $\Delta_{-1}u=\chi(D)u=\mathcal{F}^{-1}(\chi\mathcal{F}u),$
$\Delta_qu=\varphi(2^{-q}D)u=\mathcal{F}^{-1}(\varphi(2^{-q}\cdot)\mathcal{F}u),
$ for $q\geq0,$ and $S_q u=\sum^{q-1}_{i=-1}\Delta_i u=\chi(2^{-q}D)u=\mathcal{F}^{-1}(\chi(2^{-q}\cdot)\mathcal{F}u).$
\begin{definition5}\label{Def3.1}
 (Besov spaces) Assume $s\in \R,\ 1\leq p,r\leq\infty$. The nonhomogeneous Besov
 space $B^s_{p,r}(\mathbb{R}^d)$ ($B^s_{p,r}$ for short) is defined by
\begin{equation*} B^s_{p,r}=\{u\in\mathcal{S'}(\mathbb{R}^d): \|u\|_{B^s_{p,r}}=\|2^{qs}\Delta_qu\|_{l^r(L^p)}=\|(2^{qs}\|\Delta_qu\|_{L^p})_{q\geq -1}\|_{l^r}<\infty\}.\end{equation*}
In particular, $B^\infty_{p,r}=\bigcap_{s\in
\R}B^s_{p,r}.$
\end{definition5}
\begin{definition5}
Let $T>0,s\in \R$ and $1\leq p\leq \infty$. We define
\begin{equation*}
E^s_{p,r}(T):=C([0,T];B^s_{p,r})\cap C^1([0,T];B^{s-1}_{p,r}), \quad
\mbox{if} \ r<\infty,
\end{equation*}
\begin{equation*}
E^s_{p,\infty}(T):=L^\infty(0,T;B^s_{p,\infty})\cap
\mbox{Lip}([0,T];B^{s-1}_{p,\infty}),\quad
\mbox{and}\quad E^s_{p,r}:=\bigcap_{T>0}E^s_{p,r}(T).
\end{equation*}
\end{definition5}
Next we list some useful properties for Besov space $B^s_{p,r}$.
 \begin{proposition5}\label{Pro5.1}
Let $s\in \R,\ 1\leq p,r,p_i,r_i\leq\infty,i=1,2$. Then\\
(i) Topological properties: $B^s_{p,r}$ is a Banach space which is continuously embedded in $\mathcal{S'}$.\\
(ii) Density: if $1\leq p,r<\infty$, then $\mathcal{C}^\infty_c$ is
dense in $B^s_{p,r}$.\\
(iii) Embedding: if $p_1\leq p_2$ and $r_1\leq r_2$, then $B^s_{p_1,r_1}\hookrightarrow
B^{s-d(\frac{1}{p_1}-\frac{1}{p_2})}_{p_2,r_2}$. If $s_1<s_2$, then the embedding $B^{s_2}_{p,r_2}\hookrightarrow B^{s_1}_{p,r_1}$ is locally
compact.\\
(iv) Algebraic properties: if $s>0$, $B^s_{p,r}\cap L^\infty$ is a Banach
algebra. Moreover, $B^s_{p,r}$ is an algebra $\Leftrightarrow B^s_{p,r}\hookrightarrow L^\infty$
$\Leftrightarrow s>\frac{d}{p}$ or $s\geq\frac{d}{p}$ and $r=1.$\\
(v) Fatou lemma: if $\{u_n\}_{n\in \mathbb{N}}$ is a bounded sequence in
$B^s_{p,r}$, then there exist an element $u \in B^s_{p,r}$ and a subsequence $\{u_{n_k}\}_{k\in \mathbb{N}}$ such that\\
$$\lim_{k\rightarrow \infty}u_{n_k}= u \ \mbox{in} \ \mathcal{S'}, \quad \mbox{and} \quad \|u\|_{B^s_{p,r}}\leq \liminf_{k\rightarrow \infty}\|u_{n_k}\|_{B^s_{p,r}}.$$
(vi) Interpolation: (1) if $u\in B^{s_1}_{p,r}\cap B^{s_2}_{p,r}$ and $\theta\in[0,1]$, then $u\in B^{\theta s_1+(1-\theta)s_2}_{p,r}$ and
$\|u\|_{B^{\theta s_1+(1-\theta)s_2}_{p,r}}\leq
\|u\|^\theta_{B^{s_1}_{p,r}}\|u\|^{1-\theta}_{B^{s_2}_{p,r}}.$
(2) if $u\in B^{s_1}_{p,\infty} \cap B^{s_2}_{p,\infty}$ and $s_1<s_2$, then $u\in B^{\theta s_1+(1-\theta)s_2}_{p,1}$ for all
$\theta\in(0,1)$ and there exists a constant $C$ such that
$\|u\|_{B^{\theta s_1+(1-\theta)s_2}_{p,1}}\leq \frac{C}{\theta(1-\theta)(s_2-s_1)}
\|u\|^\theta_{B^{s_1}_{p,\infty}}\|u\|^{1-\theta}_{B^{s_2}_{p,\infty}}.$\\
(vii) Action of Fourier multipliers on Besov spaces: let $m\in \R$
and $f$ be a $S^m$-multiplier ($i.e.,$ $f:\R^d\rightarrow \R$ is a
smooth function and satisfies that for each multi-index $\alpha$,
there exists a constant $C_\alpha$ such that $|\partial^\alpha
f(\xi)|\leq C_\alpha(1+|\xi|)^{m-|\alpha|}$, for $\forall \xi \in
\R^d$). Then the operator $f(D)=\mathcal{F}^{-1}(f\mathcal{F})\in Op(S^{m})$ is continuous from $B^s_{p,r}$ to
$B^{s-m}_{p,r}$.\\
(viii) 1-$D$ Moser-type estimate: (1) For $s>0$, $\|uv\|_{B^s_{p,r}}\leq C (\|u\|_{B^s_{p,r}}\|v\|_{L^\infty}+\|u\|_{L^\infty}\|v\|_{B^s_{p,r}}).$
(2) If $\forall s_1\leq \frac{1}{p} <s_2$ ($s_2\geq \frac{1}{p}$ if $r=1$) and $s_1+s_2>0$, then $\|uv\|_{B^{s_1}_{p,r}}\leq C \|u\|_{B^{s_1}_{p,r}}\|v\|_{B^{s_2}_{p,r}}.$\\
(ix) The paraproduct is continuous from
$B^{-\frac{1}{p}}_{p,1}\times (B^{\frac{1}{p}}_{p,\infty}\cap
L^\infty)$ to $B^{-\frac{1}{p}}_{p,\infty}$, $i.e.$
\begin{eqnarray*}
\|uv\|_{B^{-\frac{1}{p}}_{p,\infty}}\leq
C\|u\|_{B^{-\frac{1}{p}}_{p,1}}\|v\|_{B^{\frac{1}{p}}_{p,\infty}\cap
L^\infty}.
\end{eqnarray*}
(x) A logarithmic interpolation inequality:
\begin{eqnarray*}
\|u\|_{B^{\frac{1}{p}}_{p,1}}\leq
C\|u\|_{B^{\frac{1}{p}}_{p,\infty}}\ln \big(
e+\frac{\|u\|_{B^{1+\frac{1}{p}}_{p,\infty}}}{\|u\|_{B^{\frac{1}{p}}_{p,\infty}}}\big).
\end{eqnarray*}
\end{proposition5}
Now we state the following transport equation theory that is crucial
to our purpose.
\begin{lemma5}\label{Lem5.1} (A priori estimate)
Let $1\leq p,r\leq +\infty, s>-\min\{\frac{1}{p},1-\frac{1}{p}\}$. Assume that $f_0\in B^s_{p,r}$, $F\in L^1(0,T;B^s_{p,r})$, and
$\partial_x v$ belongs to $L^1(0,T;B^{s-1}_{p,r})$ if $s>1+\frac{1}{p}$
 or to $L^1(0,T;B^{\frac{1}{p}}_{p,r}\cap L^\infty)$ otherwise. If $f\in L^\infty(0,T;B^{s}_{p,r})\cap
C([0,T];\mathcal{S'})$ solves the following transport equation
 \begin{equation}\label{Eq5.1}
\left\{\begin{array}{ll}\partial_tf+v\partial_x f=F,\\
 f\big|_{t=0}=f_0,\\
\end{array}\right.
\end{equation}
then there exists a constant $C$ depending only on $s,p,r$ such that
the following statements hold for $t\in [0,T]$\\
(i) If $r=1$ or $s\neq 1+\frac{1}{p}$,
\begin{eqnarray*}
\|f(t)\|_{B^{s}_{p,r}}\leq
\|f_0\|_{B^{s}_{p,r}}+\int_0^t\|F(\tau)\|_{B^{s}_{p,r}}d\tau
+C\int_0^tV'(\tau)\|f(\tau)\|_{B^{s}_{p,r}}d\tau,
\end{eqnarray*}
 or
 \begin{eqnarray*}
\|f(t)\|_{B^{s}_{p,r}}\leq
e^{CV(t)}(\|f_0\|_{B^{s}_{p,r}}+\int_0^te^{-CV(\tau)}
\|F(\tau)\|_{B^{s}_{p,r}}d\tau),
\end{eqnarray*}
where $V(t)=\int_0^t\|\partial_xv(\tau,\cdot)\|_{B^{\frac{1}{p}}_{p,r}\cap L^\infty}d\tau$ if $s<1+\frac{1}{p}$, and $V(t)=\int_0^t\|\partial_xv(\tau,\cdot)\|_{B^{s-1}_{p,r}}d\tau$ otherwise.\\
(ii) If $f=v$, then for all $s>0$, the estimate in (i) holds with
$V(t)=\int_0^t\|\partial_xv(\tau,\cdot)\|_{
L^\infty}d\tau.$\\
(iii) If $r<\infty$, then $f\in C([0,T];B^{s}_{p,r})$. If $r=\infty$,
then $f\in C([0,T];B^{s'}_{p,1})$ for all $s'<s$.
\end{lemma5}

\begin{lemma5}\label{Lem5.2}
(Existence and uniqueness) Let $p,r,s,f_0$ and $F$ be as in the
statement of Lemma \ref{Lem5.1}. Suppose that $v\in
L^\rho(0,T;B_{\infty,\infty}^{-M})$ for some $\rho>1,M>0$ and
$\partial_x v\in L^1(0,T;B^{\frac{1}{p}}_{p,\infty}\cap L^\infty)$
if $s<1+\frac{1}{p}$, and $\partial_x v\in L^1(0,T;B^{s-1}_{p,r})$
if $s>1+\frac{1}{p}$ or $s=1+\frac{1}{p}$ and $r=1$. Then (\ref{Eq5.1}) has a unique solution $f\in
L^\infty(0,T;B^{s}_{p,r})\cap
(\cap_{s'<s}C([0,T];B^{s'}_{p,1}))$ and the corresponding
inequalities in Lemma \ref{Lem5.1} hold true. Moreover, if
$r<\infty$, then $f\in C([0,T];B^{s}_{p,r}).$
\end{lemma5}

\begin{lemma5}\label{Lem5.3}
Let $\{v^n\}_{n\in \bar{\mathbb{N}}}, \bar{\mathbb{N}}=\mathbb{N}\cup \{\infty\}$ be a sequence of functions
belonging to $C([0,T];B^{\frac{1}{2}}_{2,1}).$ Assume that $v^n$ is
the solution to
\begin{equation*}
\left\{\begin{array}{ll}\partial_tv^n+a^n\partial_x v^n=f,\\
 v^n\big|_{t=0}=v_0,\\
\end{array}\right.
\end{equation*}
with $v_0\in B^{\frac{1}{2}}_{2,1},f\in
L^1(0,T;B^{\frac{1}{2}}_{2,1})$ and that, for some $\alpha\in
L^1(0,T),$
\begin{eqnarray*}
\sup_{n\in\mathbb{N}}\|\partial_xa^n(t)\|_{B^{\frac{1}{2}}_{2,1}}\leq
\alpha(t).
\end{eqnarray*}
If in addition $a^n$ tends to $a^\infty$ in
$L^1(0,T;B^{\frac{1}{2}}_{2,1})$, then $v^n$ tends to $v^\infty$ in
$C([0,T];B^{\frac{1}{2}}_{2,1}).$
\end{lemma5}

\begin{lemma5}\label{Lem5.4}
 (Osgood lemma) Let $\rho$ be a measurable function from $[t_0,T]$
 to $[0,a]$, $\gamma$ a locally integrable function from $[t_0,T]$
to $\R^+,$ and $\mu$ a continuous and nondecreasing function from
$(0,a]$ to $\R^+.$ Assume that, for some positive real number $c,$
the function $\rho$ satisfies
\begin{eqnarray*}
\rho(t)\leq c+\int_{t_0}^t\gamma(t')\mu(\rho(t'))dt'\quad
\mbox{for}\ \mbox{a.e.}\ t\in [t_0,T].
\end{eqnarray*}
Then we have, for a.e. $ t\in [t_0,T],$
\begin{eqnarray*}
-\mathcal{M}(\rho(t))+\mathcal{M}(c)\leq
\int^t_{t_0}\gamma(t')dt'\quad \mbox{for}\
\mathcal{M}(x)=\int_x^a\frac{dr}{\mu(r)}.
\end{eqnarray*}
\end{lemma5}

\end{document}